\RequirePackage{fix-cm}
\documentclass[smallextended]{svjour3}       
\smartqed  
\usepackage{ifpdf,appendix,caption}
\usepackage{amsmath,url}
\usepackage{amssymb}
\usepackage{enumerate}
\usepackage{graphicx}
\usepackage{color,verbatim}
\usepackage{mathrsfs,subeqnarray,subfigure}
\usepackage{listings}
\usepackage[ruled,vlined]{algorithm2e}

\newcommand{\nn}{\nonumber}

\graphicspath{{pic/}} 

\newcommand{\Rmn}[1]{\uppercase\expandafter{\romannumeral#1}}

\numberwithin{equation}{section}
\numberwithin{table}{section}
\numberwithin{figure}{section}
\newcommand{\reff}[1]{(\ref{#1})}

\newcommand{\mbR}{\mathbb{R}}
\newcommand{\be}{\begin{equation}}
\newcommand{\ee}{\end{equation}}
\newcommand{\st}{\mathrm{s.t.}}

\newcommand{\bit}{\begin{itemize}}
\newcommand{\eit}{\end{itemize}}

\newtheorem{assumption}{Assumption}
\begin{document}

\title{A Quadratic Penalty Method for    Hypergraph  Matching}

\author{Chunfeng Cui  \and
        Qingna Li \and Liqun Qi \and Hong Yan
}

\institute{Chunfeng Cui \at
              Department of Electronic Engineering, City University of Hong Kong, Kowloon, Hong Kong, China.
  This author's research was supported by the Research Grants Council (RGC)
  of Hong Kong (Project C1007-15G).
             \\
              \email{chunfcui@cityu.edu.hk}
           \and
           Qingna Li \at
           School of Mathematics and Statistics, Beijing Key Laboratory on MCAACI, Beijing Institute of Technology,
               Beijing, 100081, China.  This author's research was supported by NSFC 11671036.
             \\
              \email{qnl@bit.edu.cn}
                \and
              Liqun Qi \at
              Department of Applied Mathematics, The Hong Kong Polytechnic University,
              Hung Hom, Kowloon, Hong Kong, China.
  This author's research was supported by the Research Grants Council (RGC)
  of Hong Kong (Project C1007-15G).
  \\
              \email{maqilq@polyu.edu.hk}
              \and
              Hong Yan \at
              Department of Electronic Engineering, City University of Hong Kong, Kowloon, Hong Kong, China.
  This author's research was supported by the Research Grants Council (RGC)
  of Hong Kong (Project C1007-15G).
  \\
              \email{h.yan@cityu.edu.hk}
    }

\date{Received: 28 February 2017 / Accepted: 19 October 2017\\
@Springer Science+Business Media, LLC 2017\\
DOI 10.1007/s10898-017-0583-0}

\maketitle

\begin{abstract}
Hypergraph matching is a fundamental problem in computer vision.
 Mathematically, it maximizes a polynomial objective function,
subject to assignment constraints. In this paper, we reformulate the hypergraph matching problem as
a sparse constrained   optimization problem.
By dropping the sparse constraint, we show that the resulting relaxation problem can recover
the global minimizer  of the original problem.
 This property heavily depends on the   special structures
of  hypergraph matching.
The critical step in solving the original problem is to identify the location of nonzero entries
(referred to as the  support set) in a global minimizer. Inspired by such observation,
we  apply the quadratic penalty method to solve the relaxation problem.    Under reasonable assumptions,
we show that the support set of the global minimizer in a hypergraph matching problem can be
correctly identified   when the number of iterations is sufficiently large.
 A projected gradient method is applied as a  subsolver  to solve the quadratic penalty  subproblem.
  Numerical results demonstrate that the exact recovery of the support set indeed happens,
  and the proposed algorithm  is efficient in terms of both accuracy and CPU time.


\end{abstract}

\section{Introduction}
\label{intro}

Recently, hypergraph matching has  become a  popular tool in establishing
 correspondence between two sets of points.
It is a central problem
 in computer vision, and has been used to solve several applications,
 including
object detection \cite{berg2005shape},
image retrieval \cite{yan2005efficient},
image stitching \cite{zaragoza13,zaragoza14},
and
bioinformatics \cite{wu2012prl}.

From the point of view of graph theory,
hypergraph matching belongs to  bipartite matching.
Traditional graph matching models only use point-to-point features or
pair-to-pair features, which can be solved by
 linear assignment algorithms
\cite{jiang2007,maciel2003}
or quadratic assignment algorithms
\cite{egozi2013,Jiang2016,lee2011,litman2014,yan2015}, respectively.
To use more geometric  information
such as angles, lines, and areas,
  triple-to-triple graph matching  was proposed in 2008  \cite{zass2008probabilistic},
  and was further  studied in \cite{duchenne2011tensor,lee2011hyper,nguyen2016efficient}.
Since three vertices are associated with one edge, it
is also termed as  hypergraph matching.
Numerical   experiments in literature  \cite{duchenne2011tensor,lee2011hyper,nguyen2016efficient,zass2008probabilistic}
  show that hypergraph matching
is more efficient than   traditional graph matching.
 The aim of this paper is to study the hypergraph matching problem
    in both theory and algorithm.

The mathematical model of hypergraph matching is to
 maximize a multi-linear objective function  subject to the  row permutation constraints for $n_1\le n_2$
  \begin{equation}  \label{equ:singlestochastic}
    \begin{array}{ll}
    \Pi^1 =\{X\in\{0,1\}^{n_1\times n_2} :    \sum\limits_{l_2=1}^{n_2}X_{l_1l_2}=1, \
      l_1=1,\ldots,n_1\},
       \end{array}
  \end{equation}
  or permutation constraints for $n_1=n_2$
    \begin{equation}  \label{equ:doublestochastic}
    \begin{array}{ll}
    \Pi^2 =&\{X\in\{0,1\}^{n_1\times n_1}\, : \   \sum\limits_{l_2=1}^{n_1}X_{l_1l_2}=1, \
    l_1=1,\ldots,n_1; \\&
        \sum\limits_{l_1=1}^{n_1}X_{l_1l_2}=1, \     l_2=1,\ldots,n_1\}.
      \end{array}
  \end{equation}
  We call a matrix satisfying \reff{equ:singlestochastic} or \reff{equ:doublestochastic}
   a binary assignment matrix.
  Optimization problems   over  binary assignment matrices are known to be NP-hard due to the combinatorial property.

Most existing algorithms for hypergraph matching relax the binary constraints
into bound constraints  and solve a continuous
optimization problem. For instance,
 the probabilistic Hypergraph Matching method (HGM)     \cite{zass2008probabilistic}
  reformulated  the   constraints
 as the intersection of three convex sets,
and successively projected the variables onto the  sets until convergence.
The Tensor Matching  method (TM) \cite{duchenne2011tensor}
solved the optimization problem using the  power iteration algorithm.
The   Hypergraph Matching method  via Reweighted Random Walks (RRWHM)  \cite{lee2011hyper}
 dealt with  the problem by
 walking among two feasible vectors randomly.
Different from the above algorithms,
 Block Coordinate Ascent Graph Matching (BCAGM) \cite{nguyen2016efficient}
applied a   block coordinate ascent framework, where they kept the binary constraints, and proposed to
reformulate the multi-linear objective function into a linear one and
solve it using linear assignment algorithms.
 All the existing algorithms require the equality constraints in
\reff{equ:singlestochastic} or \reff{equ:doublestochastic} to be satisfied strictly at each iteration.
In fact,   we   only
expect  that one of the elements  is significantly larger than the others in
each  row or column of $X$. That is, the equality constraints
are only soft constraints, which allow violations to some extent.
Therefore,  we   penalize the equality
constraint violations as part of the objective function in our algorithm.

 The hypergraph matching problem can also be  reformulated equivalently as a nonlinear
 optimization problem with sparse constraint.
 During the last few years, in  the optimization community, there has been significant
 progress on solving sparse constrained nonlinear problems, particularly on dealing with
 optimality conditions and numerical algorithms in different situations.
  Recent development in optimality conditions can be found in \cite{Pan2017},
 where based on decomposition properties of the normal cones,
 the authors  characterized  different kinds of stationary points and performed detailed investigations
 on   relations of local minimizers, global minimizers and several types of stationary points.
 Other related work includes \cite{Bauschke2014Restricted,Burdakov2015MATHEMATICAL,C2016Constraint,Li2015,Pan2015On}.
The related algorithms  can be summarized into two approaches.
One is the direct approach, aiming  at dealing with the sparse constraint directly,
such as the hard-thresholding type based algorithms \cite{Beck2012Sparsity,Pan2016A}
and the $\ell_0$ penalty based algorithms \cite{Lu2012Sparse}.
 The other one is the relaxation approach such as the $\ell_p$ regularization based
 algorithms \cite{ChenLeiLuYe,Jiang2016}.
In particular, an efficient $\ell_p$ regularization algorithm was proposed
in \cite{Jiang2016}, which deals with problems over the permutation matrix
constraints \reff{equ:doublestochastic}. It can be applied to  solve
the hypergraph matching problem subject to \reff{equ:doublestochastic}.

{\bf Motivation.} Noting  that hypergraph matching  is essentially a  mixed integer programming,
 most  existing methods relax the integer constraints as box constraints, and solve the  relaxed
 continuous optimization problem. A natural question is: what is the relation between hypergraph
 matching and the relaxation problem? Furthermore, the key step in solving this problem is
 actually to identify the support set of the global minimizer. None of the existing algorithms
 has taken this fact into account. This leads to the second question: can we make use of this
 insight to design our algorithm?

{\bf Our Contributions.} In this paper, by reformulating  hypergraph matching equivalently
as a sparse constrained optimization problem, we study it from the following aspects.
\begin{itemize}

\item Relaxation problem. By dropping the sparse constraint,
we show that the relaxation problem  can recover the solution of the original problem
in the sense that the former  problem shares at least one global minimizer
 with the latter one (Theorem \ref{thm-global}).
    This result highly depends on the special structures of hypergraph matching.
    Furthermore, we show that Theorem  \ref{thm-global} can be
    extended to more general problems (Corollary 2).
    For any global minimizer of the relaxation problem, we propose a procedure to
    reduce its   sparsity until a global minimizer of the original   problem is reached.

\item Quadratic penalty method. Our aim is to identify the support set of a global
minimizer of the original problem,
 thus the equality constraints are not necessary to be  satisfied strictly.
    This motivates us to penalize the equality constraint violations, and solve
     the relaxation problem by a quadratic penalty method. We show that under
     reasonable assumptions,  the support set of a global minimizer of the original
     problem can be recovered exactly,  when the  number of iteration is
     sufficiently large (Theorems \ref{thm-penalty-original} and   \ref{thm-penalty-original-more}).

\item Projected gradient method.
For the quadratic penalty subproblem, which is a nonlinear problem with simple
 box constraints, we choose one of the active set based methods called
 the projected gradient method as a subsolver. The advantage of the active set
 based method is that it well  fits  our motivation, which is to  identify
 the support set of the solution rather than  to look for the magnitude.
 Numerical results demonstrate that the exact recovery of the support set indeed happens,
 and the proposed algorithm  is particularly suitable for
  large-scale problems.

\end{itemize}

{\bf Organization.}
The rest of the paper is organized as follows.
 In Section \ref{sec-reformulation},  we introduce the reformulation of
 the hypergraph matching problem,
  and discuss several preliminary properties.
  In Section \ref{sec-relaxation}, we study the properties of the  relaxation
  problem by dropping the sparse constraint.
  In Section \ref{sec-penalty},  we  study  the quadratic penalty method
  by penalizing the equality constraint
  violations
 and establish the convergence results in terms of support set under different situations.
An existing   projected gradient method is also discussed to solve the quadratic penalty subproblem.
Numerical experiments are reported in Section \ref{sec-numerical}.
Final conclusions are  drawn  in Section \ref{sec-conclusions}.

{\bf Notations.}
 For $x\in\mathbb{R}^n$, define the active set as $\mathcal I(x)=\{l:\,x_l=0\}$
 and the support set as $\Gamma(x)=\{l:\,x_l>0\}$.
We also use $\mathcal I^k$ and $\Gamma ^k$, and $\mathcal I^*$ and $\Gamma^*$ to
denote the corresponding sets at $x^k$ and $x^*$, respectively.
Let $|\mathcal{I}|$  be  the number of elements in  the set $\mathcal{I}$.
$\|x\|$ denotes the $\ell_2$ norm of $x$,
$\|x\|_0$  the  number of nonzero entries in   $x$,  and
$\|x\|_{\infty}$  the infinity norm of   $x$.


\section{Problem Reformulation}\label{sec-reformulation}

In this section, we will reformulate  hypergraph matching
 as a sparse constrained optimization problem, and discuss several preliminary properties.

\subsection{Hypergraph matching problem}

In this part, we will give the mathematical formulation for
 hypergraph matching, including its objective function and constraints.

Consider two hypergraphs $G_1=\{V_1, E_1\}$, and $G_2=\{V_2, E_2\}$,
where $V_1$ and $V_2$ are sets of points with $|V_1|=n_1$, $|V_2|=n_2$,
and $E_1$, $E_2$ are sets of hyperedges.
In this paper, we always suppose that $n_1\le n_2$, and
each point in $V_1$ is matched to exactly one point in $V_2$,
while each point in $V_2$ can be matched to arbitrary number of points in $V_1$.
That is, we focus on  \reff{equ:singlestochastic}.
For each hypergraph,  we consider three-uniform hyperedges. Namely,
the  three  points involved in each hyperedge are different,
for example, $(l_1,j_1,k_1)\in E_1$. Our aim is to find the best correspondence
(also referred to as `matching')   between $V_1$ and $V_2$  with the maximum matching score.

Let $X\in\mathbb{R}^{n_1\times n_2}$ be the assignment matrix between $V_1$ and $V_2$, i.e.,
\[
X_{l_1l_2}=\left\{
\begin{array}{ll}
1,& \hbox{ if }l_1\in V_1\hbox{ is assigned to }l_2\in V_2;\\
0,& \hbox{ otherwise.}
\end{array}
\right.
\]
Two hyperedges $(l_1,j_1,k_1)\in E_1$ and $(l_2,j_2,k_2)\in E_2$ are said to be matched if $l_1, j_1, k_1 \in V_1$
are assigned to $l_2, j_2, k_2\in V_2$, respectively.
 It can be  represented equivalently by
 $X_{l_1l_2}X_{j_1j_2}X_{k_1k_2}=1$.
 Let
  $\mathcal B_{l_1l_2j_1j_2k_1k_2} $ be the matching score between
  $(l_1,j_1,k_1) $ and $(l_2,j_2,k_2)$.
  Then  $\mathcal B \in \mathbb{R}^{n_1\times n_2\times n_1\times n_2\times n_1\times n_2}$  is  a
   sixth order tensor.
   Assume $\mathcal B$ is given, satisfying $\mathcal B_{l_1l_2j_1j_2k_1k_2}\ge0$
   if $(l_1,j_1,k_1)\in E_1$ and $(l_2,j_2,k_2)\in E_2$, and $\mathcal B_{l_1l_2j_1j_2k_1k_2}=0$, otherwise.

Given hypergraphs $G_1=\{V_1, E_1\}$,  $G_2=\{V_2, E_2\}$, and the matching score $\mathcal B$,
the hypergraph matching problem takes the following form
  \begin{equation}\label{equ:BX}
    \max_{X\in\Pi^1} \ \ \sum_{\tiny
    \begin{array}{c}
    (l_1,j_1,k_1)\in E_1\\
    (l_2,j_2,k_2)\in E_2
    \end{array}}\mathcal{B}_{l_1l_2j_1j_2k_1k_2} X_{l_1l_2}X_{j_1j_2}X_{k_1k_2}.
    \end{equation}
Note that  \reff{equ:BX} is a matrix optimization problem,
which can be reformulated as a vector optimization problem as follows.

Let $n = n_1n_2$, $x\in\mathbb{R}^n$ be the vectorization of $X$, that is
    \[
   x: = (\bar x_1^T,\ldots, \bar x_{n_1}^T)^T, \  \hbox{with }X = \left[
   \begin{array}{ccc}
   X_{11} & \cdots & X_{1n_2}\\
  \vdots&\ddots&\vdots\\
  X_{n_11}&\cdots& X_{n_1n_2}
   \end{array}
   \right] :=
   \left[
   \begin{array}{c }
   \bar x_1^T\\
  \vdots \\
  \bar x_{n_1}^T
   \end{array}
   \right].
    \]
Here, $\bar x_i\in \mbR^{n_2}$ is  the $i$-th block of $x$.
In the following, for  any vector $z\in \mathbb{R}^n$, we always assume it has the same partition as $x$.
 Define $\mathcal A\in\mathbb{R}^{n\times n\times n}$  as
\begin{equation}\label{equ:A}
  \mathcal{A}_{ljk} =  \mathcal{B}_{l_1l_2j_1j_2k_1k_2},
\end{equation}
where
\be\label{index-A}
  l=(l_1-1)n_2+l_2 ,
 \ j=(j_1-1)n_2+j_2 ,\
   k= (k_1-1)n_2+k_2.
\ee
Consequently,  (\ref{equ:BX}) can be reformulated as
\begin{equation}\label{equ:Ax1}
\begin{array}{rl}
\min\limits_{x\in\mathbb{R}^n} &\quad f(x):=-\frac16 \mathcal Ax^3\\
\hbox{s.t.} &\quad e^T\bar x_i = 1, \ i = 1,\ldots, n_1,\\
&\quad x\in\{0,1\},
\end{array}
\end{equation}
where  $e\in\mathbb{R}^{n_2}$ is   a vector with all entries equal to one,
and $\mathcal Ax^3:= \sum\limits_{l,j,k} \mathcal{A}_{ljk}x_lx_jx_k$.

\subsection{Preliminary properties}

In this subsection, we will discuss several  properties of $\mathcal A$, $\mathcal B$, and $f(x)$.
  We begin with properties of $\mathcal B$.

\begin{proposition}\label{prop-graph}
\begin{itemize}
\item [(i)] $\mathcal{B}_{l_1l_2j_1j_2k_1k_2}\ge0$  for all $l_1,j_1,k_1=1,\ldots,n_1$
 and $l_2,j_2,k_2=1,\ldots,n_2$;

\item [(ii)]
If $(l_1,j_1,k_1)\in E_1$, then $l_1$, $j_1$, and $k_1$ are distinct.
If $(l_2,j_2,k_2)\in E_2$, then  $l_2$, $j_2$, and $k_2$ are also distinct;

\item [(iii)] For any permutation operator $\pi$, suppose $ \pi(l_1,j_1,k_1)=(l_1',j_1',k_1')$
   and $\pi(l_2,j_2,k_2)=(l_2',j_2',k_2')$. There is
    \be\label{B-equal}
\mathcal{B}_{l_1l_2j_1j_2k_1k_2} = \mathcal{B}_{l_1'l_2'j_1'j_2'k_1'k_2'}.
\ee
\end{itemize}
\end{proposition}

The above properties of $\mathcal B$ result in the following properties of $\mathcal A$  directly.

\begin{proposition}\label{prop-A}
\begin{itemize}
     \item[(i)] $\mathcal A_{ljk}\ge0$, for all $l,j,k=1,\ldots, n$;

     \item[(ii)] For nonzero entries of $\mathcal A$, say $\mathcal A_{ljk}$,
     $x_l$, $x_j$ and $x_k$ come  from different blocks of $x$;

     \item[(iii)] Suppose $(l',j',k') $ is any permutation of $(l,j,k)$. Then
    \be \label{A-special}
    \mathcal A_{ljk} =\mathcal A_{l'j'k'}.
    \ee
\end{itemize}
        In other words,   $\mathcal{A}$  is nonnegative and symmetric.
\end{proposition}

  \textbf{Proof.}
(i)  follows directly from the nonnegativity of $\mathcal{B}$.
 In terms of (ii),  by the definition of $\mathcal A$,
there exist $(l_1,j_1,k_1)$ and $(l_2,j_2,k_2)$ such that (\ref{equ:A}) and (\ref{index-A}) hold.
Further,  we know that $x_l$ is the $l_2$-th entry in the $l_1$-th block of $x$, i.e., $x_l=(\bar x_{l_1})_{l_2}$.
Similarly,     $x_j=(\bar x_{j_1})_{j_2}$ and $x_k=(\bar x_{k_1})_{k_2}$.
By (ii) in Proposition \ref{prop-graph}, $l_1, j_1, k_1$ are distinct, which implies that $x_l$, $x_j$, $x_k$
come from different blocks of $x$.
In terms of (iii),
  since $\pi(l_1,j_1,k_1) =(l_1',j_1',k_1')$, $\pi(l_2,j_2,k_2) =(l_2',j_2',k_2')$,
  again by the definition of $\mathcal A$ and (\ref{index-A}), there is
  $B_{l_1'l_2'j_1'j_2'k_1'k_2'}=A_{l'j'k'}$.
  Together with \reff{B-equal} and (\ref{equ:A}), there is  (\ref{A-special}).
\qed

Different from other nonlinear problems,  the homogenous polynomial  $f(x)$
enjoys special structures. To see this, for  the $i$-th block $\bar x _i$, denote
\[
x _{-i} = (\bar x _1^T, \ldots, \bar x _{i-1}^T, \bar x _{i+1}^T, \ldots, \bar x _{n_1}^T)^T,\
I(i, n_2) =\{ (i-1) n_2 + 1,\ldots, i  n_2\}.
\]
Rewrite $f(x )$ as follows:
\begin{eqnarray}
\nn f(x )&=& -\frac16\sum_{l,j,k} \mathcal A_{ljk}   x_l  x_j  x_k\\
\nn&=&- \sum_{l\in I(i,n_2) \hbox{ or } j\in I(i,n_2) \hbox{ or } k\in I(i,n_2), l<j<k} \mathcal A_{ljk}  x_l  x_j x_k \\
\nn &&\quad - \sum_{l,j,k\notin I(i,n_2),l<j<k} \mathcal A_{ljk} x_l  x_j  x_k\\
&:=& f^i(\bar x_i, x_{-i}) + f^{-i}(x_{-i}). \label{equ:fi}
\end{eqnarray}

\begin{proposition}\label{prop-f0}
\begin{itemize}
  \item[(i)]
  For each block  $\bar x_i$,  $i\in\{1,\ldots,n_1\}$, $f(x)$ is a linear function of  $\bar x_i$, i.e.,
  $\nabla_{\bar x_i} f(x)$ is independent of $\bar x_i$;

  \item[(ii)]
\be\label{grad-i}
f^i(\bar x_i, x_{-i}) = \bar x_i^T\nabla_{\bar x_i} f(x).
\ee
\end{itemize}
\end{proposition}

{\bf Proof.}
In terms of (i), by the definition of $\mathcal A$, we only need to consider the term $\mathcal A_{ljk}x_lx_jx_k$,
where $\mathcal A_{ljk}$ is nonzero.
Due to (ii) in Proposition \ref{prop-A},
$\mathcal A_{ljk}x_lx_jx_k$ is linear in each related block
$\bar x_{l_1}$, $\bar x_{j_1}$, and $\bar x_{k_1}$.
Therefore, $f(x)$ is a linear function of  $\bar x_i$, $i = 1, \ldots, n_1$.

In terms of (ii),
the elements of gradient $\nabla f(x)$ take  the following form
\[
(\nabla f(x))_l = -\sum_{l<j<k}\mathcal A_{ljk}x_jx_k-\sum_{j<l<k}\mathcal A_{ljk}x_jx_k-\sum_{j<k<l}\mathcal A_{ljk}x_jx_k.
\]
Rewrite $f^i(\bar x_i, x_{-i})$ in \reff{equ:fi} as
\begin{align*}
    f^i(\bar x_i, x_{-i})= - \sum_{l\in I(i,n_2)}\left(\sum_{l<j<k}\mathcal A_{ljk}x_lx_jx_k+
    \sum_{j<l<k}\mathcal A_{ljk}x_lx_jx_k+
     \sum_{j<k<l}\mathcal A_{ljk}x_lx_jx_k\right).
\end{align*}
 Hence, $f^i(\bar x_i, x_{-i})=\sum\limits_{l\in I(i,n_2)}x_l(\nabla f(x))_l$, which gives \reff{grad-i}.   \qed

Equation (\ref{grad-i})   will be useful in Section \ref{sec-relaxation}.

\subsection{Sparse constrained optimization problem}

 Problem \reff{equ:Ax1} is a 0-1 mixed integer programming,
  which is one of Karp's 21 NP-complete problems \cite{Karp1972}.
 In this subsection, we will reformulate \reff{equ:Ax1}  into a sparse constrained optimization problem.

  By direct computations, \reff{equ:Ax1} can be reformulated as the following sparse constrained minimization problem
 \be\label{prob-hyper}\begin{array}{rl}
      \min\limits_{x\in\mathbb{R}^n}\quad  &f(x)\\
      \st\quad &  e^T\bar x_i = 1,\ i=1,\ldots, n_1,\\
      &x\ge0,\quad  \| x \|_0\le  n_1.
            \end{array}
    \ee
 To see this,  for each $x$ satisfying the  equality constraints, we have   $\| x \|_0\ge  n_1$.
 Together with $\| x \|_0\le  n_1$, we actually have $\| x \|_0= n_1$.

 In particular, if $n_1 = n_2$, by the permutation constraints \reff{equ:doublestochastic},
  problem \reff{prob-hyper} reduces to the following
   hypergraph matching problem
     \be\label{prob-1same-d}\begin{array}{rl}
      \min\limits_{x\in\mathbb{R}^n}\quad& f(x) \\
      \st\quad &  e^T\bar x_i = 1,\ i=1,\ldots, n_1,\\
      &\hat e_i^T x = 1, \ i = 1, \ldots,  n_1,\\
      & x\ge0,\quad \| x \|_0\le   n_1,      \end{array}
      \ee
  where    $\hat e_i = ((e^{n_1}_i)^T,\ldots, (e^{n_1}_i)^T)^T\in\mathbb{R}^{n},$  and
     $e^{n_1}_i$ is the $i$-th column of the $n_1$-by-$n_1$ identity matrix.

\begin{remark}
Note that the dimension of $x$ is $n=n_1n_2$,
which can be large  even for moderate  $n_1$ and $n_2$. For instance,
if $n_1=100$ and $n_2=100$, then $n=10^4$,
and the number of elements in $\mathcal{A}$ will  be around
$10^{12}$. Hence,  algorithms capable of dealing with large-scale problems are
highly in demand.
\end{remark}

\begin{remark}
  Problem (\ref{prob-hyper}) is essentially a $0$-$1$ mixed integer programming.
  Each feasible point is actually an isolated feasible point, which means that
  it is a strict local minimizer and of course is a stationary point of
  (\ref{prob-hyper}). For a theoretical verification from the optimality point of view,
   please see Theorems 1 and 3 in an  earlier version of our paper \cite{CuiLiQiYan2017}.
\end{remark}


  \section{Relaxation Problem of (\ref{prob-hyper})}\label{sec-relaxation}

  In this section, we will study the relaxation problem \reff{relax1} and its connections with the original
   problem   (\ref{prob-hyper}).

  By dropping the sparse  constraint  in (\ref{prob-hyper}), we obtain the following   problem (referred to as  the relaxation problem)
    \be\label{relax1}
        \begin{array}{rl}
             \min\limits_{x\in\mathbb{R}^n}\quad& f(x)\\
            \st\quad &  e^T\bar x_i = 1,\ i=1,\ldots,   n_1, \ \  x\ge0.
        \end{array}
    \ee

As we will show later in Theorem \ref{thm-global}, although we drop the  sparse  constraint,
the relaxation problem (\ref{relax1}) still admits a global minimizer with sparsity $n_1$ due to
the special structures of  (\ref{prob-hyper}). That is, the relaxation problem (\ref{relax1}) recovers a global minimizer of (\ref{prob-hyper}).

Let $\lambda\in\mathbb R^{n_1}, \mu\in\mathbb R^n$ be the  Lagrange multipliers corresponding to  the constraints in (\ref{relax1}).
The KKT conditions of (\ref{relax1}) are
\[
\left\{\begin{array}{ll}
\nabla_{\bar x_i} f(x)-\lambda_ie-\bar \mu_i=0,\\
   \bar x_i\ge0, \ \bar\mu_i\ge0, \ \bar x_i^T\bar \mu_i=0,  \\
  e^T\bar x_i-1=0,
\end{array}\right.
\]
which are equivalent to
\[
\left\{\begin{array}{l}
\nabla_{\bar x_i} f(x)-\lambda_ie\ge0,\ \bar x_i\ge0, \  (\nabla_{\bar x_i} f(x)-\lambda_ie)^T\bar x_i=0,\\
  e^T\bar x_i-1=0,
\end{array}\right.
\]
for all $ i = 1, \ldots, n_1$.
Define the active set and the support set for the $i$-th block $\bar x_i $ as
  \be \label{block-IT}\mathcal{I}_i (x) = \{p:\  (\bar x_i )_p =0\},\ \Gamma_i(x) =\{p:\ (\bar  x _i)_p>0\}.\ee
The KKT conditions  can be  reformulated   as
\be\label{relax1-kkt-entry}
\left\{\begin{array}{lll}
(\nabla_{\bar x_i} f(x ))_p-\lambda_i=0 ,\ &\ (\bar x _i)_p>0, \ & \ p\in\Gamma _i(x), \\
(\nabla_{\bar x_i} f(x ))_p-\lambda_i\ge0 ,\ &\ (\bar x _i)_p=0,\ &\ p\in \mathcal{I} _i(x),  \\
  e^T\bar x _i-1=0,& &
\end{array}\right.
\ee
for all $i = 1, \ldots, n_1$.
The above analysis gives the following lemma.

\begin{lemma}\label{prop-relax1}
Let $x\in\mathbb{R}^n$ be a stationary point of (\ref{relax1}),
and $\lambda \in\mathbb R^{n_1}$ be the Lagrange multiplier corresponding to the equality constraints.
For all $i=1,\ldots, n_1$, we have
\bit
\item [(i)] $
\lambda _i = \min\limits_{p\in\{1,, \ldots,   n_2\}}   (\nabla _{\bar x_i} f(x ))_p, \hbox{ and }
  (\nabla _{\bar x_i} f(x ))_p=\lambda_i $ for $p\in\Gamma _i(x)$;
\item [(ii)] $f^i(\bar x _i, x _{-i}) =  \lambda _i$.

\eit
\end{lemma}

\textbf{Proof.} (i) can be obtained directly from the KKT conditions (\ref{relax1-kkt-entry}).

  In terms of  (ii), by (\ref{grad-i}),
there is
\begin{eqnarray*}
 f^i(\bar x _i, x _{-i})&=&
 (\bar x _i)^T\nabla_{\bar x _i} f(x )\\
&=& \sum_{p\in\Gamma_i(x) }  (\bar x _i)_p(\nabla_{\bar x _i} f(x ))_p\\
&=&\sum_{p\in\Gamma_i(x)}  (\bar x _i)_p\lambda_i  \ \ (\hbox{by (i)})\\
&=&   e^T  {\bar x _i} \lambda _i  \\
&=& \lambda_i ,
\end{eqnarray*}
where the last equality is due to  $e^T\bar x_i =1$.
This completes the proof.\qed

\begin{theorem}\label{thm-global}
 There exists a global minimizer $x^{*}$ of (\ref{relax1}) such that $\|x^{*}\|_0=  n_1$.
 Furthermore, $x^{*}$ is a global minimizer of (\ref{prob-hyper}).
 \end{theorem}

{\bf Proof.}
Without loss of generality, let $y^0$ be a global minimizer of (\ref{relax1}) with $\|y^0\|_0> n_1$.
Find the first block of $y^0$, denoted as $\bar y^0_i$, such that $\|\bar y^0_i\|_0>1$.
Now we choose one index $p_0$ from $\Gamma_i(y^0):=\{p: (\bar y^0_i)_p>0\}$,
and define a new point $y^1=((\bar y^1_1)^T,\ldots, (\bar y^1_{n_1})^T)^T$ as follows:
 \[
( \bar y^1_i)_p=\left\{
 \begin{array}{ll}
 1,\ & \hbox{ if }p= p_0;\\
 0,\ &\hbox{otherwise},
  \end{array}
 \right.
  \
  \text{and}
  \quad
 \bar  y^1_{i'}=\left\{
 \begin{array}{ll}
  \bar y^1_i,\ & \hbox{ if }i'=i;\\
\bar y^0_{i'},\ &\hbox{otherwise}.
  \end{array}
 \right.
 \]
Then  $y^1$ is a feasible point for (\ref{relax1}),
 and satisfies $y^1_{-i}= y^0_{-i}$.
Furthermore,  by Proposition \ref{prop-f0}, $\nabla_{\bar x_i}f(x)$ is a function of $x_{-i}$,
there is
\be
\label{grad-equal}
\nabla_{\bar x_i} f(y^0) = \nabla_{\bar x_i} f(y^1).
\ee
Next, we will show
that $f(y^1) = f(y^0)$. Indeed,
\begin{eqnarray*}
f(y^1)-f(y^0)&=& f^i(\bar y^1_i,  y^1_{-i}) +  f^{-i}(    y^1_{-i}) -f^i(\bar y^0_i,  y^0_{-i}) - f^{-i}(    y^0_{-i})\\
&=&  f^i(\bar y^1_i,   y^0_{-i}) +  f^{-i}(   y^0_{-i}) -f^i(\bar y^0_i,   y^0_{-i}) - f^{-i}(    y^0_{-i})\\
&=&   f^i(\bar y^1_i,   y^1_{-i})  -f^i(\bar y^0_i,   y^0_{-i})\\
&=& (\bar y^1_i)^T \nabla_{\bar x_i} f(y^1)- f^i(\bar y^0_i,   y^0_{-i}) \ (\hbox{by } (\ref{grad-i}))\\
&=&(\bar y^1_i)_{p_0} (\nabla_{\bar x_i} f(y^0))_{p_0} -f^i(\bar y^0_i,   y^0_{-i}) \ (\hbox{by }(\ref{grad-equal}))\\
&=& \lambda_i  -\lambda_i \  (\hbox{by Lemma \ref{prop-relax1}})\\
&=& 0.
 \end{eqnarray*}
 This gives that $y^1$ is a feasible point with $f(y^1) = f(y^0)$.
 In other words, $y^1$ is another global minimizer of  (\ref{relax1})
  with $\|y^1\|_0<\|y^0\|_0$. If $\|y^1\|_0= n_1$, let $x^{*}:=y^1$.
  Otherwise, by repeating the above process,  we can obtain a finite sequence $y^0, y^1, \ldots, y^r$,
  which are all feasible points for (\ref{relax1})  satisfying
 \[
  \|y^r\|_0<\ldots <\|y^1\|_0<\|y^0\|_0.
 \]
 Note that there are $ n_1$ blocks in $y^0\in\mathbb R^n$. After at most $  n_1$ steps,
 the process will stop. In other words, $1\le r\le   n_1$.
  The final point $y^r$ will satisfy $\|y^r\|_0=  n_1$. One can obtain a
  global minimizer $x^{*}:=y^r$  of  (\ref{relax1}) with $n_1$ nonzero elements.

 Next, we will show that   $x^{*}$ is also a global minimizer of (\ref{prob-hyper}).
Note that the feasible region of (\ref{prob-hyper}) is a subset of the feasible region of (\ref{relax1}).
  $\|x^*\|_0=  n_1$ implies that $x^{*}$ is also a feasible point for (\ref{prob-hyper}).
Together with the fact that $f(x^{*})$ attains the global minimum of (\ref{relax1}),
  we conclude that $x^{*}$ is a global minimizer of (\ref{prob-hyper}).
\qed

Theorem \ref{thm-global}  shows that  $\|x^{*}\|_0=  n_1$ is a  necessary and sufficient condition
for a global minimizer    $x^{*}$ of (\ref{relax1})  to be a global minimizer of (\ref{prob-hyper}).
We highlight this relation  in the following corollary.

\begin{corollary}
 A global minimizer  $x^{*}$ of (\ref{relax1})   is a global minimizer of (\ref{prob-hyper})
 if and only if  $\|x^{*}\|_0=  n_1$.
\end{corollary}

  A special case of Theorem \ref{thm-global} is $|\Gamma_i(x^*)|=1$,   for each  $i\in\{1,\ldots,n_1\}$.
  Then the global minimizer  $x^{*}$ of (\ref{relax1})   is a global minimizer of (\ref{prob-hyper}).

 \begin{remark}\label{rem-true-construction}
From the proof of Theorem \ref{thm-global}, one can start from any global minimizer $y^0$
of (\ref{relax1}) to reach a point $x^*$, which is a global minimizer of both (\ref{prob-hyper}) and (\ref{relax1}).
 We only need to choose one index as the location of nonzero entry in each block $\bar y_i^0$.
  Assume $p_i$ is chosen from $\Gamma _i(y^0)$. Let $\Gamma^*_i = p_i$. This will give the
  support set in the $i$-th block, which in turn determines the global minimizer $x^*$ of (\ref{prob-hyper})  by
 \[
(\bar  x^{*}_i)_p=\left\{
 \begin{array}{ll}
 1, &\quad \hbox{ if }p = p_i,\\
 0,&\quad \hbox{otherwise,}
 \end{array}
 \right.
 \]
 for each $p\in\{1,\ldots,n_2\}$ and $i\in\{1, \ldots, n_1\}$.
One particular method to choose $p_i$ is to choose  the index with the  largest value within the block.
This is actually  the  projection of $y^0$ onto the feasible set of (\ref{prob-hyper}). Here, we summarize the process in Algorithm \ref{Alg-000}.
 \end{remark}

 \begin{algorithm}
\caption{The procedure for computing the nearest binary assignment matrix} \label{Alg-000}
\begin{itemize}
\item[Step 0.]  Given $y=(\bar y_1^T, \ldots, \bar y_{n_1}^T)^T\in\mathbb{R}^n$, a global minimizer of (\ref{relax1}).
        Let $x=0\in\mathbb{R}^n$.
\item[Step 1.] For  all  $i=1,\ldots, n_1$, find $p_i \in \arg\max_p (\bar y_i)_p$,
            and let $(\bar x_i)_{p_i}=1$.
\item[Step 2.] Output $x=(\bar x_1^T, \ldots, \bar x_{n_1}^T)^T$, which is a global minimizer of (\ref{prob-hyper}).
 \end{itemize}

\end{algorithm}

   Note that HGM   \cite{zass2008probabilistic} also solves the relaxation problem (\ref{relax1}),
   whereas  TM \cite{duchenne2011tensor}  and RRWHM \cite{lee2011hyper}  solve the relaxation
   problem with the
   permutation constraints \reff{equ:doublestochastic}.
   However, none of them analyzes the connections between the original problem and the
   relaxation problem in terms of global minimizers. On contrast, the result in Theorem
  \ref{thm-global} reveals for the first time the connections between the original
  problem (\ref{prob-hyper}) and the relaxation problem (\ref{relax1}),
  which is one of the main differences of our work from existing algorithms for hypergraph matching.

Theorem \ref{thm-global} reveals an interesting connection between the
original problem (\ref{prob-hyper}) and the relaxation problem (\ref{relax1})
in terms of global minimizers. The result  heavily relies on the property of $f(x)$ in Proposition \ref{prop-f0},
as well as the equality constraints in (\ref{prob-hyper}).
It can be   extended to the following general case.

 \begin{corollary}\label{coro:general}
  Consider
      \be\label{equ:general}
        \begin{array}{rl}
             \min\limits_{x\in\mathbb{R}^n}\quad& \hat f(x)\\
            \st\quad &  e^T\bar x_i = \alpha_i,\ i=1,\ldots,   n_1, \ \  x\ge0,
        \end{array}
    \ee
    where $\alpha_i>0$,  and $\bar x_i\in\mathbb{R}^{m_i}$ with $m_i$ being  positive integers
     satisfying $\sum_{i=1}^{n_1} m_i=n$.
    Suppose that  $\hat f(x)$ satisfies Proposition \ref{prop-f0}.
Then there exists a global minimizer $x^{*}$ of (\ref{equ:general}) such that $\|x^{*}\|_0=  n_1$.
 Furthermore, $x^{*}$ is a global minimizer of the following problem
       \[
        \begin{array}{rl}
             \min\limits_{x\in\mathbb{R}^n}\quad& \hat f(x)\\
            \st\quad &  e^T\bar x_i = \alpha_i,\ i=1,\ldots,   n_1, \ \ x\ge0,\\
                    &  \|x\|_0\le n_1.
        \end{array}
    \]
 \end{corollary}

 \section{The Quadratic Penalty Method}\label{sec-penalty}

In this section, we will consider the quadratic penalty method for
the relaxation problem (\ref{relax1}). It contains three parts.
 The first part is devoted to  motivating the quadratic penalty problem and
 its preliminary properties. The second part mainly focuses on the quadratic
 penalty method and the convergence in terms of the support set.
 In the last part, we apply an existing projected gradient method  for
  the quadratic penalty subproblem.

\subsection{The quadratic penalty problem}

 Note that  (\ref{relax1})
 is a nonlinear problem with separated simplex constraints,
 which  can be solved by
 many traditional nonlinear optimization solvers
 such as fmincon in MATLAB.
 As mentioned in Section 1,   existing algorithms for hypergraph matching
 require the equality constraints in (\ref{relax1})  to  be satisfied strictly.
 On contrast, our aim here is actually to identify the support
 set of a global minimizer of (\ref{relax1}) rather than the magnitude. Once the support set is found,
 we can follow the method in Remark \ref{rem-true-construction}
 to obtain  a global minimizer of  (\ref{prob-hyper}).
 Inspired by such observations,
 we penalize the equality constraint violations as part of the objective function.
 This is another main difference  of our method  from existing  algorithms.
 It leads us to the following quadratic penalty problem
\[
\begin{array}{rl}
\min\limits_{x\in\mathbb R^n} &\quad f(x)+\frac\sigma2\sum_{i=1}^{n_1} (e^T\bar x_i -1)^2 \\
\hbox{s.t.} &\quad  x\ge0,
\end{array}
\]
where $\sigma>0$ is a  penalty parameter.
However, this problem is not well defined in general,  since for a
fixed $\sigma$ the global minimizer  will approach infinity.
We can add an upper bound to make the feasible set bounded.
This gives the following problem
  \be\label{relax2}
    \begin{array}{rl}
        \min\limits_{x\in\mathbb{R}^n}\quad &\theta(x):=f(x)+\frac\sigma2\sum_{i=1}^{n_1} (e^T\bar x_i -1)^2\\
      \st\quad &0\le x\le M,
    \end{array}
  \ee
where {$M\ge1$} is a given  number.    (\ref{relax2})  is  actually the
quadratic penalty problem of the following problem
   \[
        \begin{array}{rl}
         \min\limits_{x\in\mathbb{R}^n}\quad&f(x)\\
        \st\quad &  e^T\bar x_i= 1,\ i=1,\ldots,   n_1,\  0\le x\le M,
        \end{array}
    \]%
    which is equivalent to (\ref{relax1}).

Having introduced the quadratic penalty problem (\ref{relax2}),  next we will
analyze the properties of   (\ref{relax2}) and its connection with the relaxation problem (\ref{relax1}).

The Lagrangian function of (\ref{relax2}) is
\[
L(x, w,\nu) = \theta(x)-  x^Tw-(M-x)^T\nu,
\]
{where $w$ and $\nu$ are the Lagrange multipliers corresponding to the inequality constraints in (\ref{relax2}).}
The KKT conditions are
\[
\left\{\begin{array}{l}
\nabla_{\bar x_i} \theta(x)-\bar w_i+\bar \nu_i=0,\\
   \bar x_i\ge0, \ \bar w_i\ge0, \ \bar x_i^T\bar w_i=0,\\
  \bar\nu_i\ge0, \  M-\bar x_i \ge0, \  \bar\nu_i^T(M-\bar x_i)=0,
\end{array}\right.
\]
for each $i\in\{1, \ldots, n_1\}$.
In particular, for a stationary point $x $ of (\ref{relax2}), let $\mathcal{I}_i(x) $ and $\Gamma_i(x)$ be defined   by (\ref{block-IT}). Define
 \[
\widehat{\Gamma} _i(x) = \{p: (\bar x _i)_p\in(0,M), p \in\Gamma _i(x)\}, \  \overline{\Gamma} _i(x) = \{p: (\bar x _i)_p=M, p \in\Gamma _i(x)\}.
\]
The KKT conditions are equivalent  to the following, for each $i\in\{1,\ldots,n_1\}$,
\be\label{quadratic-kkt}
\left\{\begin{array}{ll}
(\nabla_{\bar x_i} f(x ))_p+\sigma (e^T\bar x _i-1)\ge0,\ (\bar x _i)_p=0, & p\in \mathcal{I} _i(x), \\
 (\nabla_{\bar x_i} f(x ))_p+\sigma (e^T\bar x _i-1)=0,\ (\bar x _i)_p\in(0,M),  & p\in\widehat{\Gamma} _i(x), \\
(\nabla_{\bar x_i} f(x ))_p+\sigma (e^T\bar x _i-1)\le0,\ (\bar x _i)_p=M,  & p\in\overline{\Gamma} _i(x).
\end{array}\right.
\ee

Define the violations of the equality constraints $h\in\mathbb{R}^{n_1}$ as
\begin{equation}\label{def:h}
  h _i = e^T\bar x _i -1, \ i=1,\ldots,n_1.
\end{equation}
There  is
\be\label{h-rangle}\sigma h _i\in \left[-\max_{p\in\{1,\ldots, n_2\}} ( \nabla_{\bar x_i} f(x ))_p,\ -\min_{p\in\{1,\ldots, n_2\}} ( \nabla_{\bar x_i} f(x ))_p\right].\ee
The above analysis can be stated in the following lemma.
\begin{lemma}\label{prop-relax2}
Let $x\in\mathbb R^n$ be a stationary point of (\ref{relax2}). We have
$h _i \ge0$  for all $i = 1,\ldots, n_1$.
\end{lemma}

{\bf Proof.} For each  $i$, consider two cases.
If $\mathcal{I} _i(x)\cup{\widehat{\Gamma}} _i(x) \neq \emptyset$, by  (\ref{quadratic-kkt}),
there exists $p\in \mathcal{I} _i(x)\cup{\widehat{\Gamma}} _i(x)$  such that
$\sigma h_i  \ge -(\nabla_{\bar x_i} f(x ))_p$. By the nonnegativity of the entries in $\mathcal A$ and $x$,
there is $ -\nabla_{\bar x_i} f(x )\ge0$ and  $h_i \ge0$.
If  $\mathcal{I} _i(x)\cup{\widehat{\Gamma}} _i(x) = \emptyset$, then
$|\overline{\Gamma}_i(x)|=n_2$. In other words, $(\bar x_i )_p=M$  for all $p =1, \ldots, n_2$.
Then $h_i  = e^T\bar x _i-1 = n_2M-1\ge0$.  \qed

Let $u\in\mathbb{R}^n$  and $c\in\mathbb{R}^{n_1}$ be defined by
\[u_l =  \sum_{l<j<k}\mathcal A_{ljk} + \sum_{j<l<k}\mathcal A_{ljk} + \sum_{j<k<l} \mathcal A_{ljk},\ \ l =1,\ldots,n,
\]
and
\begin{equation}\label{equ:ci}
  c_i :=M^2\max_{p\in\{1, \ldots, n_2\}} (\bar u_i)_p,\ i=1,\ldots,n_1,
\end{equation}
where $\bar u_i$ is the $i$-th block of $u$.
It follows from  the  nonnegativity of $\mathcal{A}$ that $c\ge0$.
The following lemma describes the relation between the penalty parameter $\sigma$ and the violations of the equality constraints.

\begin{lemma}
For each stationary point $x $ of (\ref{relax2}), there is
\[
h_i  \le \frac{c_i}{\sigma},\ \forall \ i = 1, \ldots, n_1,
\]
 where $h_i$ is defined by \reff{def:h}, and $c_i$ is defined by \reff{equ:ci}.
\end{lemma}

{\bf Proof.} Note that $x \in[0,M]$. By the definition  $f(x)=-\frac16 \mathcal A x^3$,
we have $-(\nabla f(x ))_l \le M^2u_l.$  Together with (\ref{h-rangle}), there is
$\sigma h_i$$\le M^2\max_{p \in\{1, \ldots, n_2\}} (\bar u_i)_p$ $=c_i$.
The proof is complete. \qed

\begin{lemma}\label{prop2-relax2}
For each feasible point $x\in\mathbb R^n$ of (\ref{relax1}),
it is a stationary point of (\ref{relax2}) if and only if for all $i = 1, \ldots, n_1$, there is
\be\label{feasible-stationary}
(\nabla_{\bar x_i}f(x))_p = 0,  \ \forall \ p \in  \mathcal{I}_i(x)\cup \widehat\Gamma_i(x).
\ee
\end{lemma}

{\bf Proof.}
Let     $x$  be a  feasible point for (\ref{relax1}). There is $e^T\bar x_i-1=0$, $i = 1, \ldots, n_1$.
If $x$ is a stationary point of (\ref{relax1}), by the KKT conditions (\ref{quadratic-kkt}), we have
\[
(\nabla_{\bar x_i} f(x))_p+\sigma (e^T\bar x_i-1)\ge0,\ (\bar x_i)_p=0, \  p\in \mathcal{I}_i(x), \  i = 1, \ldots, n_1.
\]
Consequently,  $(\nabla_{\bar x_i} f(x))_p\ge -\sigma(e^T\bar x_i-1)=0$.
On the other hand,  $(\nabla_{\bar x_i} f(x))_p\le0$ due to the nonnegativity of
entries in  $\mathcal A$ and $x$. Therefore, $(\nabla_{\bar x_i}f(x))_p = 0$
for all $p \in  \mathcal{I}_i(x)$. For $p\in  \Gamma_i(x)$, there is $(\nabla_{\bar x_i}f(x))_p = 0.$ This gives (\ref{feasible-stationary}).

Conversely,  for  a   feasible point  $x$   for (\ref{prob-hyper}),
if (\ref{feasible-stationary}) holds,  the first two conditions in (\ref{quadratic-kkt})
hold by  $\bar x_i^Te-1=0$, $i = 1,\ldots, n_1$. For the third condition in  (\ref{quadratic-kkt}),
consider two cases. If $\overline\Gamma(x)=\emptyset,$ the result is trivial.
Otherwise, there is $(\nabla_{\bar x_i} f(x))_p\le0$ due to the nonnegativity of
entries in  $\mathcal A$ and $x$. The third condition holds automatically.
In both two cases, $x$ satisfies (\ref{quadratic-kkt}).
That is, $x$ is a stationary point of (\ref{relax2}).
 \qed

\subsection{A quadratic penalty method for (\ref{relax1})}\label{sec-algorithms}

Having investigated the properties of the quadratic penalty  problem,
we then solve (\ref{relax1}) by the traditional  quadratic penalty method,
i.e., by solving (\ref{relax2}) sequentially.
At each iteration, $x^k$ is a global minimizer of 
 the following problem
 \be \label{equ:penalty}
(P_k) \ \ \min_{0\le x\le M} \ \theta^{k}(x):=f(x)+\frac{\sigma_k}2\sum_{i=1}^{  n_1}(e^T\bar x_i-1)^2.
 \ee
 The quadratic penalty method is given in Algorithm \ref{Alg-l2}.
 \begin{algorithm}
\caption{Quadratic penalty  method for (\ref{relax1})} \label{Alg-l2}
\begin{itemize}

\item[Step 0.] Given an initial point $x^{0} \geq 0$,
          set the  parameter $ \sigma_0>0$. Let $k:=1$.
  \item[Step 1.]  Start from $x^{k-1}$ and solve ($P_k$) in \reff{equ:penalty} to obtain a global minimizer $x^{k}$.
   \item[Step 2.] If the termination rule is satisfied,
   project $x^{k}$ to {$\Pi^1$ in \reff{equ:singlestochastic} by Algorithm \ref{Alg-000}}. 
    Otherwise, choose $\sigma_{k+1} \ge \sigma_k$, $k = k+1$, and go to Step 1.
 \end{itemize}

\end{algorithm}

 The following theorem  addresses the convergence of the quadratic penalty method,
 which can be found in {classic} optimization books
such as   \cite[Theorem 17.1]{NumericalOpt} and  \cite[Corollary 10.2.6]{SunYuan}.
Therefore, the proof is omitted.

  \begin{theorem}\label{thm-penalty}
 Let $\{x^k\}$ be generated by Algorithm \ref{Alg-l2}, and $\lim_{k\to \infty}\sigma_k=+\infty$.
Then any accumulation point of the generated sequence $\{x^k\}$  is a global minimizer of  (\ref{relax1}).
\end{theorem}

Due to Theorem \ref{thm-penalty}, in   following analysis, we always assume the following holds.

\begin{assumption}\label{ass1}
  Let $\{x^k\}$ be generated by Algorithm  \ref{Alg-l2} and  $\lim_{k\rightarrow \infty}\sigma_k=+\infty$.
  Denote $K$ as a subset of $\{1,2,\ldots\}$.
  Assume that $\lim_{k\rightarrow \infty, k\in K}x^k=z$,
 and $z$ is a global minimizer of  (\ref{relax1}).
\end{assumption}

The next theorem mainly addresses the relation between the support set of $x^k$
and that of the global minimizer of (\ref{prob-hyper}).
Recall that for $x^k$, there is
\[
\mathcal{I}^k=\{l: \ x_l^k=0\}, \ \ \Gamma^k=\{l: \ x_l^k> 0\}.
\]

\begin{theorem}\label{thm-penalty-original}
Suppose that Assumption  \ref{ass1} holds.
If there exists a positive integer  $k_0$,  such that $\|x^k\|_0= n_1$ for all $k\ge k_0, \ k\in K$,
then  there is a positive integer $k_1\ge k_0$ such that the support set of $z$ can be  identified correctly. That is,
\[
\Gamma^k = \Gamma (z),\ \text{for all} \ k\ge k_1, \ k\in K.
\]
Furthermore, $z$ is a global minimizer of (\ref{prob-hyper}).

\end{theorem}

{\bf Proof.}  First, we show $|\Gamma(z)|=  n_1$. Noting that $z$ is a global minimizer of (\ref{relax2}), we have
 \[
 \Gamma(z)\ge  n_1.
 \] Since $\lim_{k\to+\infty,k\in K }x^k= z$, there exists a positive integer $k'$ such that  for $k\ge k', \ k\in K$, there is
  \[
 x^k_l>1/2, \text{ for all } l\in\Gamma(z).
 \]
 This implies that $\Gamma(z)\subseteq \Gamma^k$.
 It follows from the assumption that $|\Gamma^k|=  n_1$ for all $k\ge k_0, \ k\in K $.
   Consequently, we have  $|\Gamma(z)|=  n_1$.   Therefore,     $\Gamma(z)=\Gamma^k$ holds for $k\ge k_1:=\max\{k_0,k'\}$.
The second part  holds following the second part of Theorem \ref{thm-global}.
 The proof is finished.\qed

 Theorem \ref{thm-penalty-original} indicates that we do not need to drive $\sigma_k$ to infinity since
 only  the support set of $z$ is  needed. If the conditions in  Theorem \ref{thm-penalty-original}
 hold, then we can stop the algorithm when the number of elements in  $\mathcal{I} ^k$ keeps unchanged for
 several iterations.
However, if  there is  $\|x^k\|_0>n_1$, we need more notations to analyze  the connections.

Let $\mathcal{J}^k_i$ be the set of indices  corresponding to the largest values in
the $i$-th block $\bar x^k_i$, $p^k_i$ be the smallest index in $\mathcal{J}^k_i$, and $\mathcal{J}^k$ be the set of indices
containing the largest values in each block of $x^k$, i.e.,
\[
\mathcal{J}^k_i = \arg\max_{p} \{(\bar x^k_i)_p\}, \  p^k_i = \min \{p: p\in \mathcal{J}^k_i\}, \hbox{ and } \mathcal{J}^k:=\bigcup_{i = 1}^{n_1} \{p_i^k + n_2(i-1)\}.
\]
Similarly, we define
\[
\mathcal{J} _i (z)= \arg\max_{p} \{(\bar z_i)_p\}, \  p _i(z) = \min \{p: p\in \mathcal{J} _i(z)\},  \ \mathcal{J}(z):=\bigcup_{i = 1}^{n_1} \{p_i (z) + n_2(i-1)\}.
\]

\begin{theorem} \label{thm-penalty-original-more}
Suppose that  Assumption \ref{ass1}  holds.
\bit
\item [(i)] If $\|z\|_0= n_1$,
 then there exists an   integer $k_0>0$, such that $\Gamma(z) = \mathcal{J}^k$  for all $k\ge k_0, \ k\in K$;
\item [(ii)] If $\|z\|_0>n_1$ and $ |\mathcal{J} _i(z)|=1$ for all $i = 1, \ldots, n_1$,
then there exists  a global minimizer $x^{*}$ of (\ref{prob-hyper}) and a positive integer $k_0$,
such that for all $k\ge k_0, k\in K$, there is $\Gamma^{*} = \mathcal{J}^k$;
\item [(iii)] If $\|z\|_0>n_1$ and $ |\mathcal{J} _i(z)|> 1$ for one $i = 1, \ldots, n_1$,
 then there exists  a global minimizer $x^{*}$ of (\ref{prob-hyper}),
 a subsequence $\{x^k\}_{k\in K'}$ and a positive integer $k_0$,
 such that for all $k\ge k_0, k \in K'$, there is $\Gamma^{*} = \mathcal{J}^k$.
\eit
\end{theorem}

{\bf Proof.}  With Theorem \ref{thm-global} and $\|z\|_0=n_1$, $z$ must be a global minimizer of (\ref{prob-hyper}).
By the definition of $\Gamma(z)$ and $\mathcal I(z)$,  there exists an  integer
$k_0>0$, such that for all $k\ge k_0, k\in K$, there is
$z_l >z_{l'}$ for ${l'}\in \mathcal{I}(z)$ and $l \in \Gamma(z)$.
 This gives $\mathcal{J}^k_i = \mathcal{J} _i(z)$ and (i).

  In terms of  (ii),  $|\mathcal{J} _i(z)|=1$ implies that for $k\in K$ sufficiently large, there is
\[
(\bar x^k_i)_{\mathcal{J}_i(z)}>(\bar x^k_i)_p, \ \forall\, p \notin\mathcal{J}_i(z),\ \  i = 1, \ldots, n_1.
\]
Consequently,  there is $\mathcal{J}^k_i= \mathcal{J}_i(z) $. Now let
 $y^0:=z$. Similar to the arguments in the proof of Theorem \ref{thm-global},
 we construct $y^1$ by choosing $p_0=\mathcal{J} _i(z)$.
 Then we can obtain a finite sequence $y^0, y^1, \ldots, y^r$  with
 \[
  \|y^r\|_0<\ldots<\|y^1\|_0<\|y^0\|_0.
 \]
After at most $  n_1$ steps, the process will stop. In other words, $1\le r\le   n_1$.
At the final point $y^r$ will satisfy $\|y^r\|_0=  n_1$.
One can find a global minimizer $x^{*}:=y^r$  of problem (\ref{relax1}) with sparsity $  n_1$.
Further, $x^{*}$ is also  a global minimizer of (\ref{prob-hyper})  and   satisfies
 \[
  |\mathcal{J}^{*}_i| =1,\ \mathcal{J}_i^{*} =\Gamma_i^*= \mathcal{J}_i(z) = \mathcal{J}_i^k.
 \]
Consequently, (ii) holds.

 For (iii),  suppose  there exists an index $q_1$ such that $|\mathcal{J} _{q_1}(z)|>1$.
Consequently, there exists $p_1\in \mathcal{J}^k_{q_1}$, such that for $k\in K$ sufficiently large,
there are infinite number of $k$ satisfying $\mathcal{J}_{q_1}^k = p_1$.
Denote the corresponding subsequence as $\{x^k\}_{k\in K_1}$, where $K_1\subset K$.
Similarly, for  $|\mathcal{J} _{q_2}(z)|>1$, we can find an infinite number of $k\in K_2\subseteq K_1$
such that $\mathcal{J}_{q_2}^k = p_2$. Repeating the process until for all  blocks,
there exists an  integer $k_0>0$, such that $|\mathcal{J}_i^k|=1$,  $i = 1,\ldots, n_1$,
for all $k\in K_t\subseteq K_{t-1}\ldots\subseteq K_1$, $ k\ge k_0$.
Let $K':=K_t$. Now similar to  Remark \ref{rem-true-construction}, for all $i = 1, \ldots, n_1$, we define $x^{*}$ as follows:
  \[
(\bar  x^{*}_i)_{p_i}=\left\{
 \begin{array}{ll}
 1, & \hbox{ if }{p_i} = \mathcal{J}_i^k,\ k\in K', \ k\ge k_0,\\
 0,&\hbox{  otherwise.}
 \end{array}
 \right.
 \]
Then we find a global minimizer of (\ref{relax1}) such that $\|x^{*}\|_{0}=n_1$.
 For $k\ge k_0$, $k\in K'$, there is $\mathcal{J}_i^k =\mathcal{J}_i^{*}$, $i = 1, \ldots, n_1$.
  Consequently, $x^{*}$ is also a global minimizer of (\ref{prob-hyper}). Hence, (iii) holds.
This completes the proof. \qed

Theorems \ref{thm-penalty-original} and \ref{thm-penalty-original-more}
state that  there is always a subsequence of  $x^k$ whose support set
 will coincide with
the support set of one global minimizer of  (\ref{prob-hyper}). Consequently, it provides
a method to design the termination rule for Algorithm \ref{Alg-l2}.

\subsection{A projected gradient method for the subproblem (\ref{relax2})}\label{sec-subsolver}

In this subsection, we will use a   projected gradient method  to solve the subproblem.

Note that the subproblem (\ref{relax2}) is  a nonlinear problem with simple  box constraints.
Various methods can be chosen to solve (\ref{relax2}), one of which is the active set based method.
We  prefer such type of method  because it quite  fits  our motivation
to identify the support set of the global minimizer of (\ref{prob-hyper}) rather than the magnitude.
The strategy of identifying the active set is therefore crucial in solving (\ref{relax2}).
We choose a popular  approach proposed in \cite{Bertsekas1982Bertsekas},
and modify it into the resulting projected gradient method, as shown in Algorithm \ref{Alg-Newton}.
Other typical projected gradient methods in \cite{calamai1987,dai2006} can also be used.

\begin{remark}
\label{rmk:supportset}
Note that the projected  gradient method is only guaranteed
to converge to a stationary point. Based on
Lemma \ref{prop-relax2},   the sum of each block in the stationary point  is larger than or equal to one.
In other words,  at least one entry in each block  is larger
 than zero. This will partly explain the  numerical observation that the magnitudes
 of the returned solution by our algorithm clearly fall  into two parts: the estimated active   part,
 which is close to zero, and the estimated nonzero part. The latter part is actually
 the estimated support set where the true support set of global minimizers of (\ref{prob-hyper}) lies in.
 Moreover, based on   Remark \ref{rem-true-construction},
 one could  identify the support set of a global minimizer of (\ref{prob-hyper}) easily.
 On the other hand, noting that the quadratic penalty  problem (\ref{relax2}) is in general nonconvex,
 it is usually not easy to find a global minimizer. Fortunately, our numerical results demonstrate that in many cases,
 the projected gradient method  can return  a solution with  accurate support set.
\end{remark}

Note that the relaxation problem {(\ref{relax1})} does not take any sparsity into account.
However, as shown in Theorem  \ref{thm-global}
at least one of the global minimizers of the relaxation problem  {(\ref{relax1})}
is a  global minimizer  of the original problem  {(\ref{prob-hyper})}.
By  the quadratic penalty method, we can indeed identify the support set of
one global minimizer of  {(\ref{prob-hyper}) }under reasonable assumptions.

\begin{remark}
We focused on the problem (\ref{prob-hyper}) so far. One may wonder whether the
theoretical results can be extended to  (\ref{prob-1same-d}).
It turns out that the extension is not   trivial and the analysis becomes more
challenging and complicated due to the equality constraints $e^T\bar x_i - 1 = 0$ and $\hat e_i^Tx -1 = 0$.
We leave it as a  topic to study in future.
 However, as we will demonstrate in the numerical part, the algorithm designed here
 can also be applied to solving  the relaxation problem of (\ref{prob-1same-d}).

\end{remark}

\begin{algorithm}
\caption{Projected  gradient  method}\label{Alg-Newton}
\begin{itemize}

\item[Step 0.] Given an initial point $x^0\in \mathbb{R}^{n}$ with $0 \leq x^0  \leq M$
               and tolerance  $\hbox{Tol} >0$.
               Set the parameters as $0<\rho<1/2$, $\epsilon > 0$,
               $0< \beta <1$,  $M\ge 1$.
               Let $j:=0.$ Denote   $P(x)$ as the projection of $x\in\mathbb R^n$ onto the box constraint $0\le x\le M$, and $g(x) = \nabla \theta(x)$.
       \item[Step 1.] Calculate the estimated active set at $x^j$ as
    \[
 I_j:= \left\{l\;|\;0\leq x_l^j\leq\epsilon^j, \;g_l(x^j)>0;   \ \hbox{or } M-\epsilon^j\le  x^j_l\le M, \; g_l(x^j)<0;\;l=1,\ldots, n  \right\},
\]
 where $\epsilon ^j = \min\{\epsilon,\omega ^j\}$,  $\omega ^j =
\|x^j - P(x^j- N g(x^j))\|$,
and $N$ is a fixed positive
definite diagonal  matrix in $\mathbb{R}^{n\times n}$. Let $\bar{I}_j := \{1,\ldots, n\}
 \backslash I_j$.

  \item[Step 2.] Calculate the residual $\delta^j\in \mathbb{R}^{n}$ by
  \[
\delta^j: =\left[\begin{array}{l}
\delta^j_{ I_j}  \\
\delta^j_{\bar I_j}
\end{array}\right]
 \]
  with $\delta^j_{ I_j} = \min\left\{x^j_{ I_j},\;g_{ { I_j}} (x^j)\right\}$ and
  $\delta^j_{\bar I_j} = g_{ {\bar I_j}}(x^j) $. If $\|\delta^j\|\leq \hbox{Tol}$, stop. Otherwise, go to Step
3.

 \item[Step 3.] Calculate the direction $d^j\in \mathbb{R}^{n}$ by
  \begin{equation*}
d^j:= \left[\begin{array}{l}
d^j_{ I_j}  \\
d^j_{\bar I_j}
\end{array}\right] \end{equation*}
  where $d^j_{ I_j} = -(Z^j)^{-1} x_{ I_j}$,
    $Z^j \in \mathbb{R}^{|I_j|\times |I_j|}$ is a   positive definite diagonal matrix,
    and $d^j_{\bar I_j}=-\eta^jg_{{ I_j}}(x^j)$, where $\eta^j$ is a scaling parameter.

  \item[Step 4.] Choose the step size as $\alpha^j = \beta^{m^j}$, where $m^j$ is
  the smallest nonnegative integer  $m$  such that the following condition holds
 \[
  \theta(P(x^j+\beta^{m}d^j))- \theta(x^j) \leq \rho \left(\beta^{m}\sum_{l\in \bar I_j}g_l(x^j)d^j_l+ \sum_{l\in I_j}g_l(x^j)(P(x^j_l+\beta^{m}d^j_l)-x_l^j)\right).
  \]
  \item[Step 5.] Update $x^{j+1}$ by $x^{j+1}=P(x^j+\alpha^jd^j)$, 
  $j:=j+1$. Go to Step 1.
\end{itemize}

\end{algorithm}

 \section{Numerical Results}\label{sec-numerical}

In this section, we will evaluate the performance of our algorithm and compare it  with
 several state-of-the-art   approaches for  hypergraph matching.

\subsection{Implementation issues}
\label{subsec:implement}
Our algorithm is termed  as  QPPG, which is the abbreviation of
 Quadratic Penalty Projected Gradient method. Basically, we run Algorithm \ref{Alg-l2}
 (referred to as outer iterations) and solve the subproblem \reff{equ:penalty}
 by calling Algorithm \ref{Alg-Newton} (referred to as inner iterations).
   In practice, we only execute
 an inexact version of Algorithm 2 by one step.
   QPPG2  means that Algorithm \ref{Alg-l2} is applied  to permutation constraints \reff{equ:doublestochastic}.
  For TM \cite{duchenne2011tensor},
 RRWHM \cite{lee2011hyper},
HGM \cite{zass2008probabilistic},
and BCAGM \cite{nguyen2016efficient},
we use the authors' MATLAB codes and C++  mex files.
 Our algorithm is  implemented in MATLAB (R2015a),
while tensor vector multiplications are computed with  C++ mex files.
All the experiments are preformed on a Dell desktop with Intel dual core i7-4770 CPU at
3.40 GHz and 8GB of memory running  Windows 7.

 In Algorithm \ref{Alg-l2},
 set  $\sigma_0=10$ and the initial point  $x^0$ as the  vector with all entries equal to one.
  Update  $\sigma_k$  as
\begin{equation}\label{equ:sigma}
  \sigma_{k+1}=\left\{
             \begin{array}{ll}
               \min(10^5, 1.3\sigma_k), & \hbox{ if } \sum_i |h_i^k|\ge0.1; \\
               \min(10^5, 1.2\sigma_k), & \hbox{ if } \underline h^k \le \sum_i |h_i^k|<0.1; \\
                \sigma_k, & \hbox{ otherwise,}
             \end{array}
           \right.
\end{equation}
where $h_i^k=e^T\bar x_i^k-1$ and  $\underline h^k$ is the maximal value of $\sum_i|h_i^k|$ for five consecutive steps.
We stop Algorithm  \ref{Alg-l2} if one of the following conditions is satisfied:
(a) $|\Gamma^k|$ is less than $1.2 n_1$;
(b) $|\Gamma^k|$  stays    unchanged for  ten consecutive steps.
  As for the output, each  $x^k$ returned by different algorithms is projected to its nearest binary assignment matrix
   by Algorithm \ref{Alg-000}   except  HGM and BCAGM, which output a binary assignment matrix directly).
  The parameters in Algorithm  \ref{Alg-Newton}  are
  $\hbox{Tol} =10^{-5}$,   $\rho=10^{-6}$,  $\epsilon =10^{-2}$,  and $\beta=0.5$.
  $\eta^j$ is chosen as $\eta^j=\frac{n_1}{\|g_{\bar{I}_j}(x^j)\|_{\infty}}$.
  The positive definite diagonal matrices $N$ and $Z^j$ are set to be the identity matrix.

{\bf Generating Tensor $\mathcal A$.}
Note that $\mathcal{A}\in \mathbb{R}^{n\times n\times n}$   contains   $n^3$ elements.
Fortunately, in hypergraph matching, as  analyzed in Proposition
\ref{prop-A}, $\mathcal A$ has   special structures.
Further, $\mathcal A$ is also sparse.
 There are three steps to generate $\mathcal A$.
The first step is to   construct   hyperedges $E_1$ and $E_2$, where
each hyperedge connects three different points.
The hyperedges in $E_1$ are generated  by randomly selecting three points in $V_1$.
 We fix $|E_1|$ as $n$.
 $E_2$ contains the nearest triples to elements in $E_1$,
 and is generated  following  the
nearest neighbour  query approach  in \cite{duchenne2011tensor,nguyen2016efficient}.
The second step is to generate $\mathcal B$.
Note that the number of nonzero entries in $\mathcal B$ are at most
$|E_1||E_2|$, which will be   large even for moderate $|E_1|$ or $|E_2|$.
In fact, for each hyperedge in $E_1$, we only use  $s$   nearest hyperedges in $E_2$  to construct $\mathcal B$.
In other words, $\mathcal B$ is calculated by
 \be\label{equ:B}
   \mathcal{B}_{l_1l_2j_1j_2k_1k_2} =\left\{
                        \begin{array}{ll}
                          \exp\{-\gamma\|f_{l_1j_1k_1}-f_{l_2j_2k_2}\|\}, \ &
                          \hbox{ if }(l_2,j_2,k_2)\in E_2   \hbox{ is one of the }s\\& \hbox{ nearest neighbours of } (l_1, j_1, k_1), \\
                          0, \ & \hbox{ otherwise,}
                        \end{array}
                      \right.
 \ee
    where $f_{l_1j_1k_1}$ and $f_{l_2j_2k_2}$ are  feature vectors determined
    by hyperedges $(l_1, j_1, k_1)$ and $(l_2, j_2, k_2)$,
     and $\gamma = \frac{1}{\text{mean}(\|f_{l_1j_1k_1}-f_{l_2j_2k_2}\|)}$      \footnote{$\hbox{mean}(\|f_{l_1j_1k_1}-f_{l_2j_2k_2}\|)=\frac{\sum_{\mathcal B_{l_1l_2j_1j_2k_1k_2}>0}\|f_{l_1j_1k_1}-f_{l_2j_2k_2}\|}{\hbox{number of } B_{l_1l_2j_1j_2k_1k_2}>0}$.} is  a normalization parameter.
    Here, for each $(l_1, j_1, k_1)\in E_1$,
the $s$ nearest neighbours are
the $s$ smallest solutions for
$\min_{(l_2, j_2, k_2)\in E_2}\|f_{l_1j_1k_1}-f_{l_2j_2k_2}\|$.
Then $\mathcal A$ can be obtained according to (\ref{equ:A}).
The number of nonzero elements is
$O(sn)$,       which is linear in   $n$.
     Therefore,  $\mathcal{A}$  is a  sparse  tensor.

We evaluate the numerical performance  mainly from the following three aspects:
  (1) `Accuracy':  denoting the ratio of successful matching, calculated by
  \[\frac{\hbox{number of correctly identified support indices }}{\hbox{number of true support indices}};\]
(2) `Matching Score':  calculated by $\frac16\mathcal A(x_B^k)^3$,
where $x_B^k$ is the nearest binary assignment  vector  of $x^k$ generated by Algorithm \ref{Alg-000};
 (3) `Running Time': the total CPU time   in seconds.
 For each algorithm (except BCAGM), we only count the computing time for solving the optimization problem.
 However,  BCAGM has to  compute all elements in  $\mathcal A$
 to
obtain   results with high accuracy.
Therefore, the running time for BCAGM contains two parts: generating $\mathcal A$  with all elements and
solving the optimization problem.

\begin{figure}[h]
  \centering

     \includegraphics[width=1.52in]{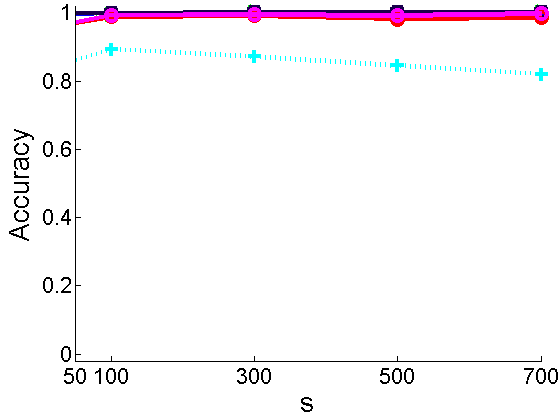}
  \includegraphics[width=1.52in]{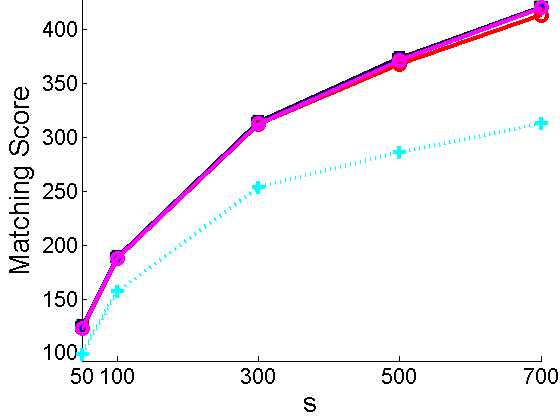}
  \includegraphics[width=1.52in]{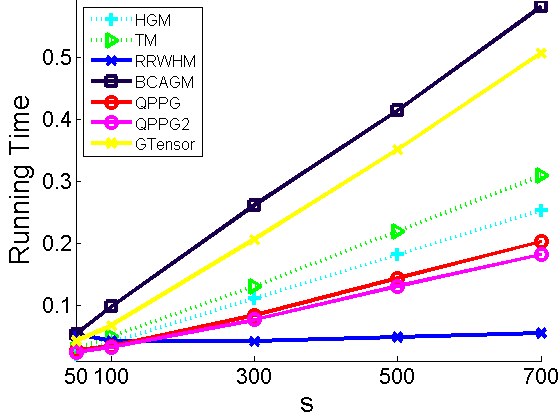}\\
   \caption{Results for  different $s$.} \label{fig:sparse}%

\end{figure}

{\bf Role of  Sparsity  of $\mathcal A$.} To see this, we test different values of $s$
on the   examples from the CMU house dataset\footnote{Downloaded from
 \text{http://vasc.ri.cmu.edu/idb/html/motion/house/}},
 which   has been    widely used in literature \cite{duchenne2011tensor,lee2011hyper,nguyen2016efficient,zhou2015local}.
 For all  examples, there is   $n_1=30$ and $n_2 = 30$.
 We take all 111 pictures with labels from 0 to  110,
 which are  the same house taken from slightly different viewpoints.
 That is, two houses with close labels are similar.
  For each picture with label  $v_1$,
 we match  it with   $v_1+60$.
 In other words, matching  picture $v_1$ with   $v_1+60$ is a test problem.
 Then we change
 $v_1$ from $0$ to $50$ to produce $51$ test examples.
 To save time of generating input data $\mathcal A$,  only elements with $l\le j \le k$ are computed in $\mathcal{A}$, and the time  consumed  is denoted by   `GTensor'.  The average results for the test examples  are reported in Figure \ref{fig:sparse}.    One can see that
 CPU time for generating tensor  is not neglectble  comparing with CPU  time for solving the problem.
 On the other hand, the accuracy stays almost unchanged for $s\ge100$.
  Note that the matching score will be larger when $s$ increases.
  It is reasonable as a denser   $\mathcal{A}$ will result in
a larger objective function.
 Therefore, we set  $s=100$ in all the following tests.

\begin{figure}[h]
  \centering
     \includegraphics[width=1.52in]{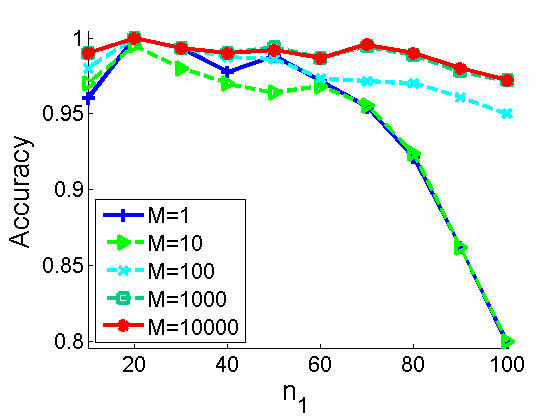}
  \includegraphics[width=1.52in]{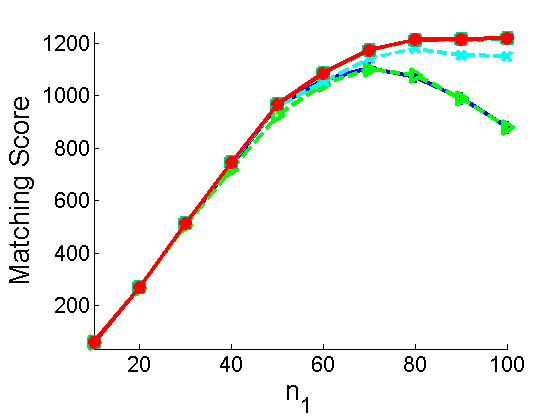}
  \includegraphics[width=1.52in]{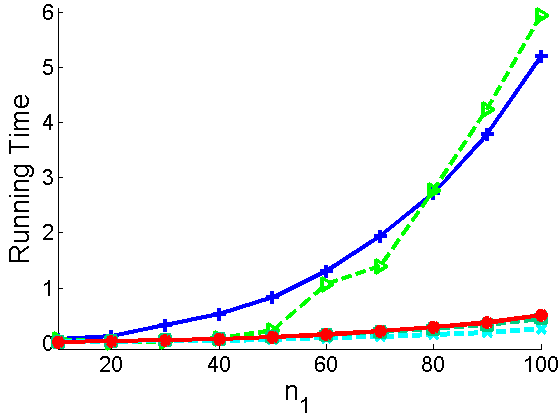}\\
   \caption{Results of QPPG for  different $M$.} \label{fig:QPPGM}%

\end{figure}

{\bf Role of Upper Bound $M$.}
 To see the role of $M$, numerical tests  are performed  on the synthetic data
 following the approach in \cite{duchenne2011tensor,nguyen2016efficient}.
Firstly, $n_1$ points in $V_1$ are sampled following the 
 standard normal  distribution
  $\mathcal{N}(0,1)$.
  Secondly,  points in   $V_2$ are computed by $V_2 = TV_1+\epsilon$, where
  $T\in\mathbb{R}^{n_1\times n_1}$ is a transformation matrix,
  and $\epsilon\in\mathbb{R}^{n_1}$ is the Gaussian noise.
 We choose $n_1=n_2$ ranging from 20 to 100, and $M$ from 1 to 10000.
 All experiments are executed for  100 times, and the  average results     are reported in Figure \ref{fig:QPPGM}.

  We can see that    $M=1000$ or $M=10000$ produces competitive results,
  while $M\le100$   is not good for large  problems in terms of both accuracy and CPU time.
     A possible reason is that   small $M$ might lead to   less flexibility
    for  the entries in $x$.
     Hence, in  the following results, we choose $M=10000$.

\subsection{Performance of QPPG and QPPG2}
\label{subsec:QPPG}

\begin{figure}[h]
  \centering
  \includegraphics[width=1.52in]{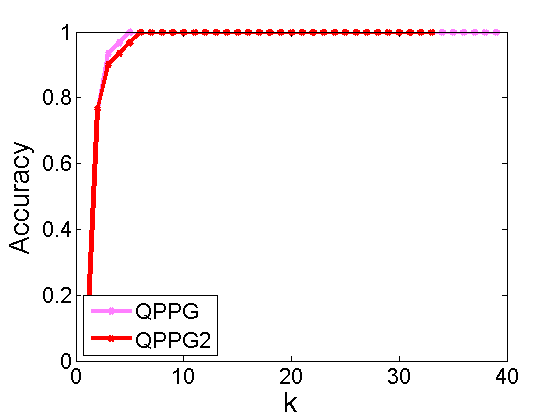}
  \includegraphics[width=1.52in]{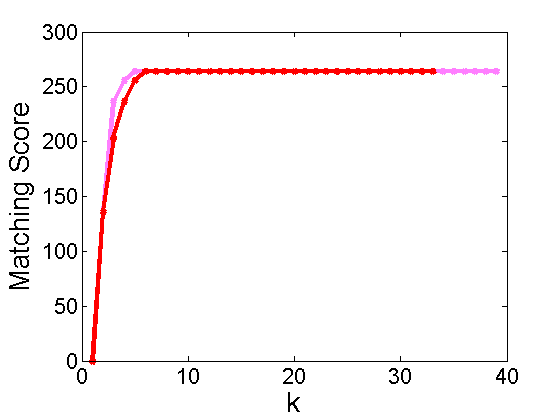}
  \includegraphics[width=1.52in]{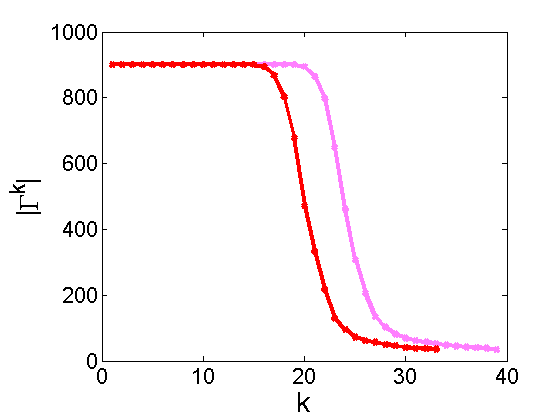}
    \caption{Accuracy, matching score and $|\Gamma^k|$
   while running QPPG and QPPG2.}
    \label{fig:stepbystep}
\end{figure}

In this subsection, we will illustrate the performance of our algorithm with
synthetic data   discussed above.
   We set   $n_1=n_2=30$.
Figure \ref{fig:stepbystep}  shows the information while running Algorithm \ref{Alg-l2},
 including  the accuracy, matching score and  size of support set
   at $x^k$.

From Figure \ref{fig:stepbystep}, one can find that $|\Gamma^k|$
keeps unchanged in the first few steps, and then drops
rapidly from $n_1^2$ to $n_1$,
while both   accuracy and matching score reach their maximum value within five
 steps. It shows the potential of our algorithm  for identifying the exact
 support set quickly, even during the process of iteration.
 This motivates us to stop our algorithm when $|\Gamma^k|$ is
 small enough, or stay unchanged for several iterations.

\begin{figure}[h]
  \centering
  \includegraphics[width=4in]{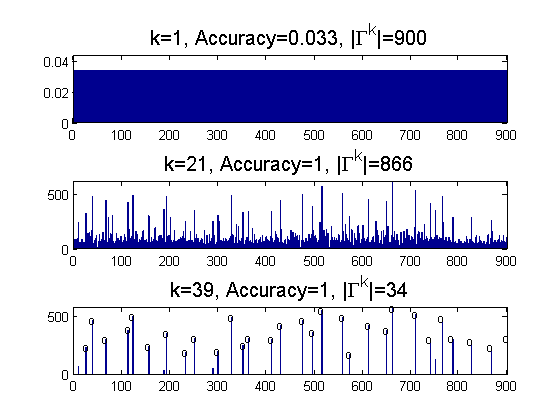}
    \caption{Entries in $x^k$ with $k = 1, 21, 39$ by QPPG.
    The small circles in the bottom figure denote the true support set.}
    \label{fig:xk_steps}
\end{figure}

We also report the magnitude of entries in  $x^k$ at several selected steps of QPPG in  Figure \ref{fig:xk_steps}.
The algorithm stops at $k=39$.  One can see that
$|\Gamma^k|$ is decreasing. At the final step,
the solution is  sparse.
This  coincides with Remark \ref{rmk:supportset}, i.e.,
 the magnitudes of the returned solution by our algorithm clearly fall into two parts:
 the estimated active   part, which is usually close to zero,
 and the estimated nonzero part, which is the support set we are looking for.

 \subsection{CMU house dataset}

\begin{figure}[h]
  \centering
     \includegraphics[width=1.52in]{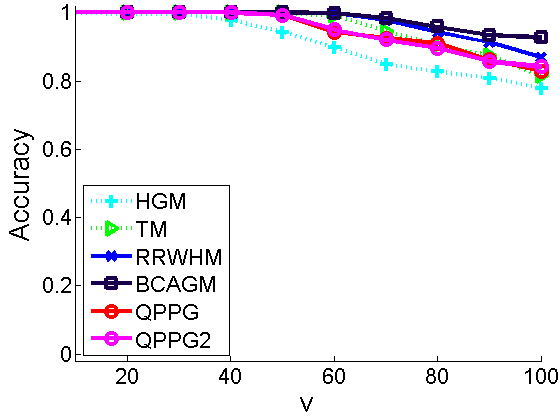}
  \includegraphics[width=1.52in]{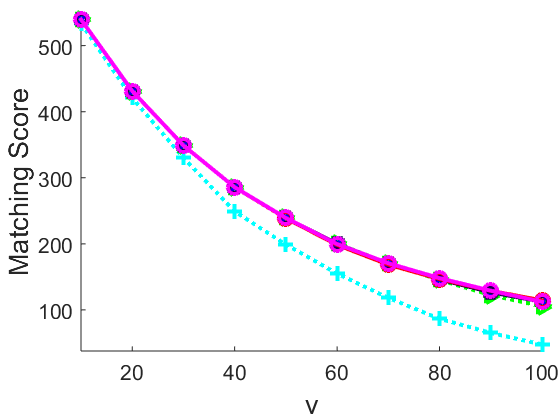}
  \includegraphics[width=1.52in]{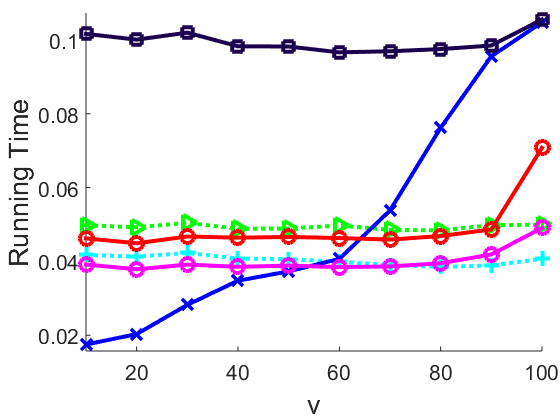}\\
   \caption{Results for CMU house dataset  with $n_1=n_2=30$.} \label{fig:30-30}%

\end{figure}

In this subsection, we will test our algorithms on
 the CMU house dataset.
 Similar to   Section \ref{subsec:implement},
we try to match picture $v_1$ with  $v_2$.  As $v_1$ and $v_2$ change,
we deal with   different hypergraph matching test problems.
 For a fixed value $v=|v_1-v_2|$,     we set $v_1=0,\ldots,110-v$ and $v_2=v,\ldots,110$.
 The total  number of test examples is $111-v$. We test these examples,
 and plot the average results for each $v$ in Figure \ref{fig:30-30}.
One can see that most algorithms (except HGM) achieve good performance
in terms of both accuracy and matching score.
In terms of CPU time,    QPPG and QPPG2 are competitive  with HGM and TM, and
  faster than other methods.

\begin{figure}[h]
  \centering
     \includegraphics[width=1.5in]{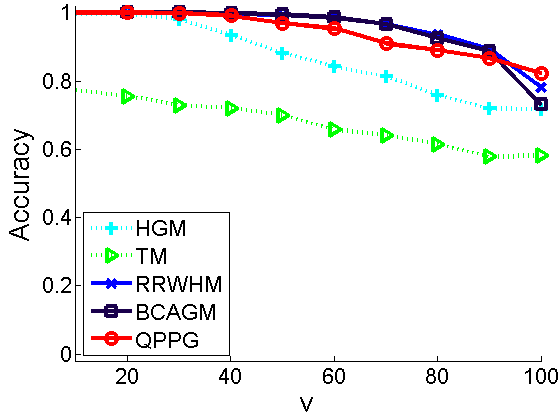}
  \includegraphics[width=1.5in]{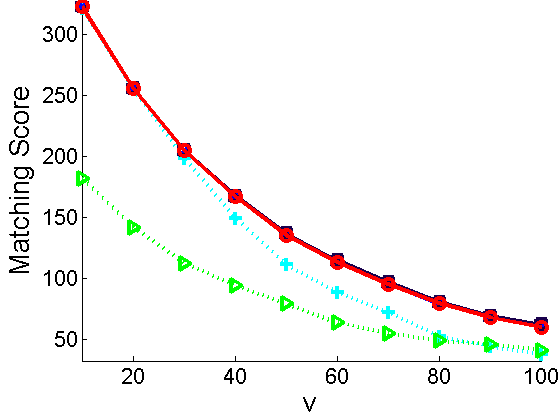}
  \includegraphics[width=1.5in]{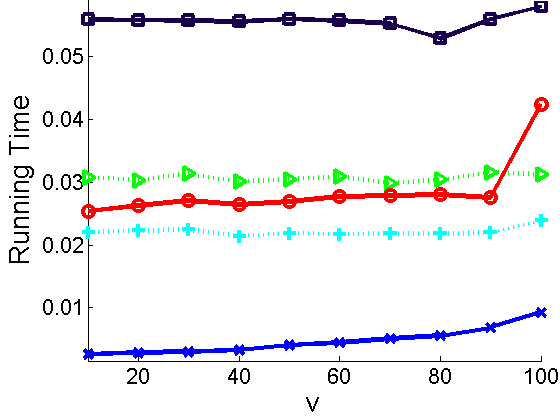}\\
   \caption{Results for CMU house dataset with $n_1=20$ and $n_2=30$.} \label{fig:20-30}%

\end{figure}

\begin{figure}[h]
  \centering
     \includegraphics[width=4in]{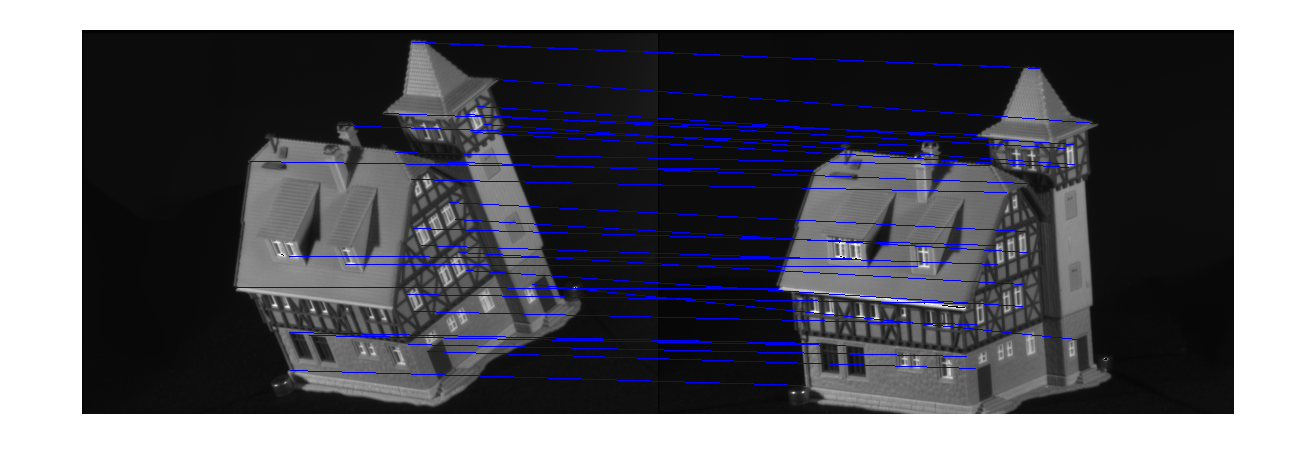}\\
  \caption{The  matching results for two houses with $v_1=0$ and $v_2=60$  by QPPG.  The blue lines are   point-to-point correspondence.
}\label{fig:matches}
\end{figure}

We also compare QPPG  with other algorithms on CMU house dataset with $n_1=20$ and $n_2=30$.
The results are obtained in a similar  way  as that for   Figure \ref{fig:30-30},
and are shown in Figure \ref{fig:20-30}.
One can see that  QPPG performs   well in both
   accuracy and matching score.
   As for CPU time, all the algorithms are competitive since   the maximum time   is about 0.06s.
 Figure \ref{fig:matches} shows the matching results for two houses
  with $v_1=0$ and $v_2=60$.

\subsection{Large dimensional synthetic data}

\begin{figure}[h]
  \centering
  \includegraphics[width=1.5in]{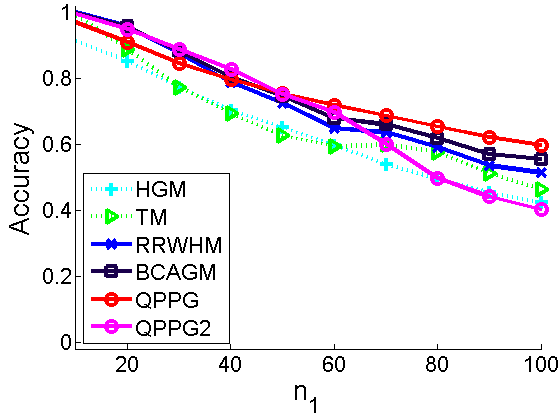}
  \includegraphics[width=1.5in]{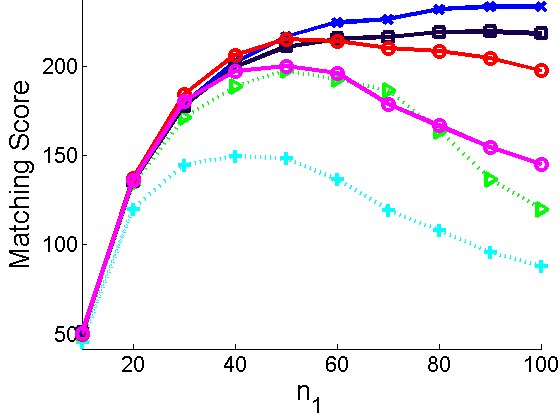}
  \includegraphics[width=1.5in]{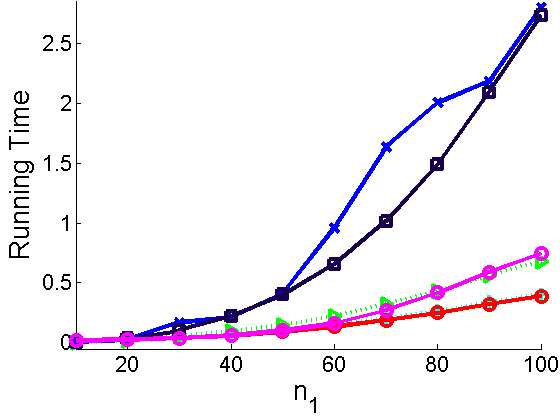}  \\
  \caption{Results for the   fish dataset.}\label{fig:res_fish}
\end{figure}

\begin{figure}[h]
  \centering
  \includegraphics[width=2.5in]{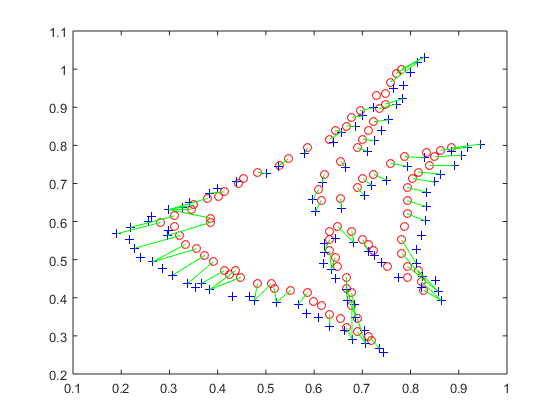}\\
  \caption{The  matching results  for fish dataset by QPPG.
  The red circles `$\circ$' stand  for   points  in  $V_1$,
and blue plus signs `$+$' represent points in $V_2$.
The green lines are   point-to-point correspondence.}\label{fig:matching_fish}
\end{figure}

In this section, large dimensional problems
in the fish dataset\footnote{Downloaded from
 http://www.umiacs.umd.edu/$\sim$zhengyf/PointMatching.htm}
 are used to test our algorithms.
We use all 100 examples in the subfolder \texttt{res\_fish\_def\_1}.
For each example, $V_1$ is the set of target fish, and $V_2$ is the set of deformation fish.
The number of points in each set is around 100.
 Our task is  to  match the two sets.
 We select $n_1=n_2=10, 20, \ldots, 100$ points randomly  from each fish
(for   fish with  less than 100 points, we use all the points).
The average results are shown in Figure \ref{fig:res_fish}.
It can be seen that
our algorithm is competitive with other methods  in terms of accuracy,
matching score and CPU time.
One of the matching  results is shown in Figure \ref{fig:matching_fish}.

\begin{figure}[h]
  \centering
  \includegraphics[width=1.5in]{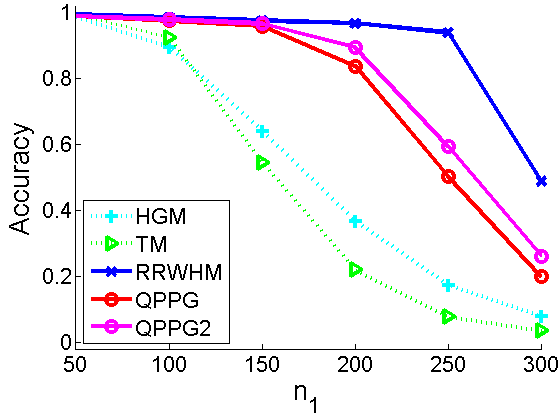}
  \includegraphics[width=1.5in]{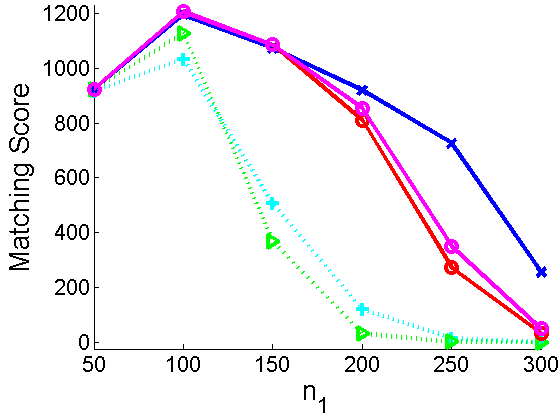}
  \includegraphics[width=1.5in]{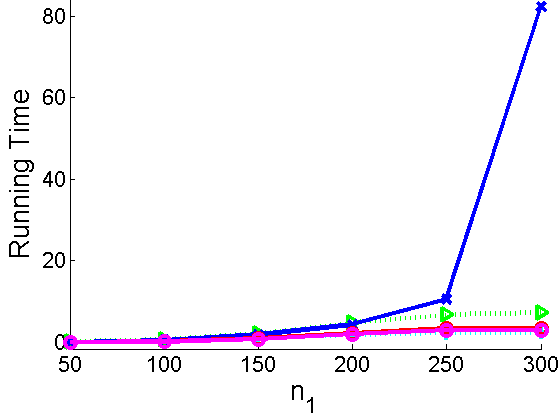}  \\
  \caption{Results for  synthetic data with large  $n_1$.}\label{fig:res_LargeDim}
\end{figure}

 Furthermore, synthetic data explained in Section  \ref{subsec:implement}
 is also used to test these algorithms.
  All the algorithms are tested except  BCAGM, as their codes run into memory troubles for large-scale problems.
 We choose $n_1=n_2$ from 50 to 300, and repeat the tests for 100 times.
  The average  results are reported in Figure \ref{fig:res_LargeDim}.
   One can see that  QPPG and QPPG2  perform comparably well with RRWHM
   in terms of both accuracy and matching score for $n_1$ less than or equal to 200.
 For $n_1$ greater than or equal to $250$,
  the running  time for QPPG and QPPG2 increases slowly  as $n_1$ increases,
 which implies that the proposed algorithm can deal with large-scale
  problems while returning   good matching results.

 \section{Conclusions}\label{sec-conclusions}

In this paper, we reformulated   hypergraph matching   as
a sparse constrained   optimization problem.
By dropping the sparse constraint,
we  showed that the relaxation problem
has at least one global minimizer, which is also the global minimizer of the
original problem.
Aiming at seeking for the support set of the global minimizer of the original problem,
we allowed   violations of the equality constraints by penalizing them in a quadratic form.
Then a quadratic penalty method
 was  applied to solve the relaxation problem. Under reasonable assumptions,
 we showed that the support set of the global minimizer in hypergraph matching can be
 identified correctly without driving the penalty parameter to infinity.
 Numerical results demonstrated the high accuracy of the support set returned by our method.

\section*{Acknowledgements}
The authors would like to thank Dr. Yafeng Liu from Academy of Mathematics and Systems Science,
Dr. Bo Jiang  from Nanjing Normal University,
and
Dr. Lili Pan from Beijing Jiaotong University for  discussions and insightful comments on this paper.
We are also grateful to two anonymous reviewers for their valuable comments,
which further improved the quality of this paper.


\begin{thebibliography}{10}

\bibitem{Bauschke2014Restricted}
{ Bauschke, H.~H., Luke, D.~R., Phan, H.~M., and Wang, X.:}
\newblock Restricted normal cones and sparsity optimization with affine
  constraints.
\newblock {Found. Comput. Math. } {\bf 14}(1), 63--83 (2014)

\bibitem{Beck2012Sparsity}
{ Beck, A., and Eldar, Y.~C.:}
\newblock Sparsity constrained nonlinear optimization: Optimality conditions
  and algorithms.
\newblock {SIAM J. Optim.} {\bf 23}(3),  1480--1509 (2013)

\bibitem{berg2005shape}
{ Berg, A.~C., Berg, T.~L., and Malik, J.:}
\newblock Shape matching and object recognition using low distortion
  correspondences.
\newblock  {IEEE Conf.  Computer Vision and Pattern Recognition}
 {\bf 1}, 26--33 (2005)

\bibitem{Bertsekas1982Bertsekas}
{ Bertsekas, D.~P.:}
\newblock Projected newton methods for optimization problems with simple
  constraints.
\newblock {SIAM J. Control  Optim.}
 {\bf 20}(2),  221--246 (1982)

\bibitem{Burdakov2015MATHEMATICAL}
{ Burdakov, O.~P., Kanzow, C., and Schwartz, A.:}
\newblock Mathematical programs with cardinality constraints: Reformulation by
  complementarity-type conditions and a regularization method.
\newblock {SIAM J.  Optim.}
{\bf 26}(1), 397--425 (2016)

\bibitem{calamai1987}
{ Calamai, P.~H., and Mor{\'e}, J.~J.:}
\newblock Projected gradient methods for linearly constrained problems.
\newblock {Math. Program.}
{\bf 39}(1), 93--116 (1987)

\bibitem{C2016Constraint}
{ {C}ervinka, M., Kanzow, C., and Schwartz, A.:}
\newblock Constraint qualifications and optimality conditions for optimization
  problems with cardinality constraints.
\newblock {Math. Program.}
{\bf 160}(1),   353--377 (2016)

\bibitem{ChenLeiLuYe}
{ Chen, X., Guo, L., Lu, Z., and Ye, J.~J.:}
\newblock An augmented lagrangian method for non-lipschitz nonconvex
  programming.
\newblock {SIAM J. Numer. Anal.}
{\bf 55}, 168--193 (2017)



\bibitem{CuiLiQiYan2017}
Cui, C. F., Li, Q. N., Qi, L. Q. and Yan, H.:
A quadratic penalty method for hypergraph matching.
arXiv:1704.04581v1 (2017)

\bibitem{dai2006}
{ Dai, Y.-H., and Fletcher, R.:}
\newblock New algorithms for singly linearly constrained quadratic programs
  subject to lower and upper bounds.
\newblock {Math. Program.}
{\bf 106}(3),  403--421 (2006)

\bibitem{duchenne2011tensor}
{ Duchenne, O., Bach, F., Kweon, I.-S., and Ponce, J.:}
\newblock A tensor-based algorithm for high-order graph matching.
\newblock {IEEE Trans. Pattern Anal. Mach. Intell.}
 {\bf  33}(12),  2383--2395 (2011)

\bibitem{egozi2013}
{ Egozi, A., Keller, Y., and Guterman, H.:}
\newblock A probabilistic approach to spectral graph matching.
\newblock {IEEE Trans. Pattern Anal. Mach. Intell.}
{\bf 35}(1), 18--27 (2013)

\bibitem{Jiang2016}
{ Jiang, B., Liu, Y.~F., and Wen, Z.:}
\newblock {$L_p$}-norm regularization algorithms for optimization over
  permutation matrices.
\newblock {SIAM J. Optim.}
{\bf 26}(4), 2284--2313 (2016)

\bibitem{jiang2007}
{ Jiang, H., Drew, M.~S., and Li, Z.-N.:}
\newblock Matching by linear programming and successive convexification.
\newblock {IEEE Trans. Pattern Anal. Mach. Intell.}
{\bf 29}(6),  959--975 (2007)

\bibitem{Karp1972}
Karp, Richard M.:
Reducibility among combinatorial problems.
Complexity of computer computations.
springer US, 85--103 (1972)

\bibitem{lee2011hyper}
{ Lee, J., Cho, M., and Lee, K.~M.:}
\newblock Hyper-graph matching via reweighted random walks.
\newblock {IEEE Conf. Computer Vision and Pattern Recognition.},
{\bf}  1633--1640 (2011)

\bibitem{lee2011}
{ Lee, J.-H., and Won, C.-H.:}
\newblock Topology preserving relaxation labelling for nonrigid point matching.
\newblock {IEEE Trans. Pattern Anal. Mach. Intell.}
{\bf 33}(2), 427--432 (2011)

\bibitem{Li2015}
{ Li, X., and Song, W.:}
\newblock The first-order necessary conditions for sparsity constrained
  optimization.
\newblock {J.  Oper. Res.  Soc.  China}
{\bf 3}(4),  521--535 (2015)


\bibitem{litman2014}
{ Litman, R., and Bronstein, A.~M.:}
\newblock Learning spectral descriptors for deformable shape correspondence.
\newblock {IEEE Trans. Pattern Anal. Mach. Intell.}
{\bf 36}(1), 171--180 (2014)

\bibitem{Lu2012Sparse}
{ Lu, Z., and Zhang, Y.:}
\newblock Sparse approximation via penalty decomposition methods.
\newblock {SIAM J.  Optim.}
{\bf 23}(4),  2448--2478 (2013)

\bibitem{maciel2003}
{ Maciel, J., and Costeira, J.~P.:}
\newblock A global solution to sparse correspondence problems.
\newblock {IEEE Trans. Pattern Anal. Mach. Intell.}
{\bf  25}(2),  187--199 (2003)

\bibitem{nguyen2016efficient}
{ Nguyen, Q., Tudisco, F., Gautier, A., and Hein, M.:}
\newblock An efficient multilinear optimization framework for hypergraph
  matching.
\newblock {IEEE Trans. Pattern Anal. Mach. Intell.}
{\bf 39}(6), 1054-1075 (2017)

\bibitem{NumericalOpt}
Nocedal, J.,  and Wright, S.:
Numerical optimization.
Springer Science \& Business Media (2006)


\bibitem{Pan2017}
{ Pan, L., Xiu, N., and Fan, J.:}
\newblock Optimality conditions for sparse nonlinear programming.
\newblock {Sci. China Math.}
  {\bf 60}(5), 759--776 (2017)

\bibitem{Pan2015On}
{ Pan, L., Xiu, N., and Zhou, S.:}
\newblock On solutions of sparsity constrained optimization.
\newblock {J.  Oper. Res.  Society of China}
{\bf 3}(4),  421--439 (2015)

\bibitem{Pan2016A}
{ Pan, L., Zhou, S., Xiu, N., and Qi, H.:}
\newblock A convergent iterative hard thresholding for sparsity and
  nonnegativity constrained optimization.
\newblock {Pac. J. Optim.}
{\bf 13}(2), 325-353 (2017)



\bibitem{SunYuan}
Sun, W. Y.,  and Yuan, Y.-X.:
Optimization theory and methods: nonlinear programming (Vol. 1).
Springer Science \& Business Media   (2006)



\bibitem{wu2012prl}
{ Wu, M.-Y., Dai, D.-Q., and Yan, H.:}
\newblock Prl-dock: Protein-ligand docking based on hydrogen bond matching and
  probabilistic relaxation labeling.
\newblock {Proteins. Struct.  Funct.  Genet.}
{\bf 80}(9), 2137--2153 (2012)

\bibitem{yan2005efficient}
{Yan, H.:}
\newblock Efficient matching and retrieval of gene expression time series data
  based on spectral information.
\newblock  {Intern. Conf.  Comput.   Sci.  Appl.},
{\bf}  357--373 (2005)

\bibitem{yan2015}
{Yan, J., Zhang, C., Zha, H., Liu, W., Yang, X., and Chu, S.~M.:}
\newblock Discrete hyper-graph matching.
\newblock {IEEE Conf.  Computer Vision and   Pattern Recognition},
{\bf\/}  1520--1528 (2015)

\bibitem{zaragoza13}
{ Zaragoza, J., Chin, T.-J., Brown, M.~S., and Suter, D.:}
\newblock As-projective-as-possible image stitching with moving DLT.
\newblock{IEEE Conf.  Computer Vision and   Pattern Recognition},
 {\bf\/}  2339--2346 (2013)

\bibitem{zaragoza14}
{ Zaragoza, J., Chin, T.-J., Tran, Q.-H., Brown, M.~S., and Suter, D.:}
\newblock As-projective-as-possible image stitching with moving DLT.
\newblock {IEEE Trans.  Pattern Anal. Mach. Intell.},
{\bf  36}(7),  1285--1298 (2014)

\bibitem{zass2008probabilistic}
{ Zass, R., and Shashua, A.:}
\newblock Probabilistic graph and hypergraph matching.
\newblock{IEEE Conf.  Computer Vision and   Pattern Recognition},
{\bf\/} 1--8 (2008)

\bibitem{zhou2015local}
{ Zhou, J., Yan, H., and Zhu, Y.:}
\newblock Local topology preserved tensor models for graph matching.
\newblock  {IEEE Conf. Syst. Man. Cybern.},
  {\bf\/}  2153--2157 (2015)

\end{thebibliography}
\end{document}